\documentclass[12pt]{article}
\usepackage{latexsym}
\usepackage{amsmath,amsthm,amssymb}
\usepackage{psfrag, graphicx}

\newtheorem{theorem}{Theorem}
\newtheorem{pr}{Proposition}[section]

\newtheorem{corollary}{Corollary}

\theoremstyle{remark} 
\theoremstyle{definition}
\newtheorem*{proof1}{Proof of Theorem \ref{thm:circinv}}
\newtheorem*{proof2}{Proof of Corollary \ref{2cor}}
\newtheorem*{proof3}{Proof of Theorem \ref{thm:spherinv}}
\newtheorem*{proof4}{Proof of Theorem \ref{thm:wavefinite}}
\newtheorem*{proof5}{Proof of \req{invP} in  Theorem \ref{thm:invPW} for $n=2$}
\newtheorem*{proof6}{Proof of Theorem \ref{thm:e2m}}
\newtheorem*{proof7}{Proof of Theorem \ref{thm:invPW} for $n>2$}

\newcommand{\eproof}{{\mbox{\ }~\hfill
\mbox{\large $\Box$} \par \vskip 10pt}}
\newcommand{\supp}{\mbox{\rm supp\,}}
\newcommand{\eps}{\varepsilon}
\newcommand{\R}{\mathbf{R}}
\newcommand{\N}{\mathbf{N}}

\DeclareMathOperator{\diam}{diam}

\DeclareMathOperator{\M}{\mathcal{M}}
\DeclareMathOperator{\W}{\mathcal{W}}
\DeclareMathOperator{\B}{\mathcal{B}}

\DeclareMathOperator{\Io}{\mathcal{I}}
\DeclareMathOperator{\Po}{\mathcal{P}}
\DeclareMathOperator{\Do}{\mathcal{D}}

\DeclareMathOperator{\Bd}{\mathbf{B}}
\DeclareMathOperator{\Id}{\mathbf{I} }
\DeclareMathOperator{\Dd}{ \mathbf{D}}
\DeclareMathOperator{\Sd}{ \mathbf{S}}

\newcommand{\ps}{p}
\newcommand{\rM}{r}
\newcommand{\rH}{\bar r}
\newcommand{\rT}{t}
\newcommand{\rP}{\rho}
\newcommand{\ro}{\rM}
\newcommand{\tio}{\rT}
\newcommand{\Rball}{R_0}
\newcommand{\Lap}{\Delta}

\newcommand{\ntw}{(n-2)/2}

\newcommand{\ph}{\varphi}

\newcommand{\tm}{ \subseteq}
\newcommand{\di}{ \partial}
\newcommand{\req}[1]{(\ref{eq:#1})}
\newcommand{\beq}[1]{\begin{equation}\label{eq:#1}\begin{aligned}}
\newcommand{\eeq}{\end{aligned}\end{equation}}
\newcommand{\beqn}{\begin{equation*}\begin{aligned}}
\newcommand{\eeqn}{\end{aligned}\end{equation*}}

\newcommand\set[1]{\left\{  #1 \right\} }
\newcommand\abs[1]{\vert  #1 \vert }
\newcommand\wabs[1]{\left\vert  #1 \right\vert }
\newcommand\f{\mathbf}

\usepackage{algorithmicx}
\usepackage[ruled]{algorithm}
\usepackage{algpseudocode}
\newenvironment{alginc}[1][pseudocode]{\medskip\algsetlanguage{#1}\begin{algorithmic}[1]}{\end{algorithmic}\medskip}

\begin{document}

\title{Inversion of spherical means and the wave equation in even dimensions}
\author{David Finch
\\Department of Mathematics
\\Oregon State University
\\Corvallis, OR 97331
\\ finch@math.oregonstate.edu
\and Markus Haltmeier
\\Department of Computer Science
\\Universit\"{a}t Innsbruck
\\Technikerstra\ss e 21a 
\\A-6020 Innsbruck, Austria
\\Markus.Haltmeier@uibk.ac.at
\and
Rakesh
\\Department of Mathematical Sciences
\\University of Delaware
\\Newark, DE 19716
\\rakesh@math.udel.edu}
\date{January 15, 2007} \maketitle
 
\begin{abstract}
  We establish inversion formulas of the so called filtered
  back-projection type to recover a function supported in the ball in
  even dimensions from its spherical means over spheres centered on
  the boundary of the ball. We also find several formulas to recover
  initial data of the from $(f,0)$ (or $(0,g)$) for the free space wave
  equation in even dimensions from the trace of the solution 
  on the boundary of the ball, provided the initial data has support in the
  ball. 
\end{abstract}

\section{Introduction and Statement of Results}

The problem of determining a function from a subset of its spherical means
has a rich history in pure and applied mathematics. Our interest in the
subject was provoked by the new medical imaging technologies called
thermoacoustic and photoacoustic tomography. The idea behind these
\cite{KruMilReyKisReiKru00, WanPanKuXieStoWan03} is to illuminate an object
by a short burst of radiofrequency or optical energy which causes rapid
(though small in magnitude) thermal expansion which generates an acoustic
wave.  The acoustic wave can be measured on the periphery or in the exterior
of the object. The inverse problem we consider is to find the distribution
of the absorbed energy throughout the body. This is of interest, since the
amount of energy absorbed at different points may be diagnostic of disease
or indicative of uptake of probes tagged to metabolic processes or gene
expression \cite{KruKisReiKruMil03}. For a more thorough discussion of the
modelling and biomedical applications, the reader is referred to the recent
survey \cite{XuMWan06}.  If the illuminating energy is impulsive in time,
the propagation may be modelled as an initial value problem for the wave
equation.  The problem of recovering the initial data of a solution of the
wave equation from the value of the solution on the boundary of a domain is
of mathematical interest in every dimension, but for the application to
thermo-/photoacoustic tomography it would appear that the three dimensional
case is the only one of interest, since sound propagation is not confined to
a lower dimensional submanifold. However, there exist methods of measuring
the generated wave field which do not rely on point measurements of the sort
that would be generated by an (idealized) acoustic transducer. In
particular, integrating line detectors, which have been studied in
\cite{BurHofPalHalSch05,PalNusHalBur06}, in effect compute the integral of
the acoustic wave field along a specified line. In this paper, we work under
the assumption that the speed of sound, $c$, is constant throughout the
body, and since the x-ray transform in a given direction of a solution of
the three dimensional wave equation is a solution of the two dimensional
wave equation, the problem is transformed. If a circular array of line
detectors is rotated around an axis orthogonal to the direction of the line
detectors \cite{HalFid06, PalNusHalBur06}, then for each fixed rotation
angle the measurement provides the trace of the solution of the two
dimensional wave equation on the circle corresponding to the array. The
initial data of this two dimensional problem is the x-ray transform of the
three dimensional initial data. If the inital data can be recovered in the
disk bounded by the detector array and assuming that the projection of the
object to be imaged lies in this disk, then the problem of recovering the
three dimensional initial data is reduced to the inversion of the x-ray
transform in each plane orthogonal to the axis of rotation. One such two
dimensional problem is illustrated in Figure \ref{fg:tact-line}.

\begin{psfrags}
\psfrag{e3}{$e_1$}
\psfrag{L}{line detector}
\psfrag{x}{$p\in S$}
\psfrag{E}{$\R^2$}
\psfrag{sigma}{rotation}
\begin{figure}[htb!]
    \begin{center}
    \includegraphics[width = 0.6\textwidth]{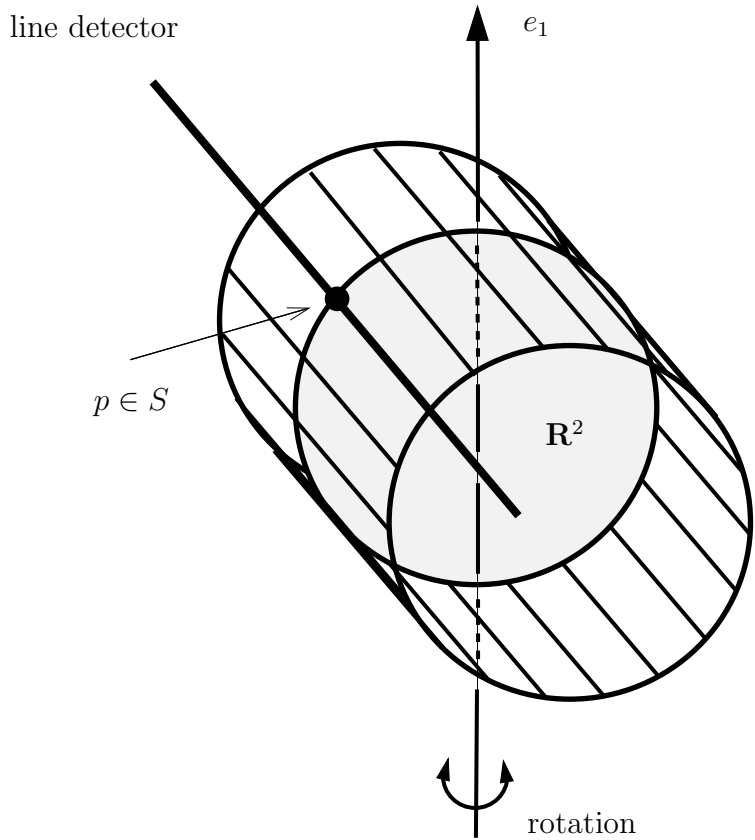}
   \end{center}
    \caption{ {\bf Principle of thermoacoustic tomography with integrating
             line detectors.}  A cylindrical array of line detectors records
             the acoustic field and is rotated around the axis $e_1$. For
             fixed rotation angle the array outputs the x-ray transform
             (projection along straight lines) of the solution of the wave
             equation restricted to lines passing through the boundary $S$
             of the disk.  The initial condition is given by the x-ray
             transform of the initially induced pressure restricted to lines
             orthogonal to the base of the cylinder.}
\label{fg:tact-line}
\end{figure}
\end{psfrags}

To our knowledge, the first work to tackle the problem of recovering a
function from its circular means with centers on a circle  was
\cite{Nor80}, whose author was interested in ultrasound reflectivity
tomography. He found an inversion method based on harmonic
decomposition and for each harmonic, the inversion of a  Hankel
transform. This method has been the basis for most subsequent work on exact
inversion of circular means. The inversion of the Hankel transform involves
a quotient of a Hankel transform of a harmonic component of the data and a
Bessel function. That this quotient be well-defined turns out to be a
condition on 
the range of the circular mean transform \cite{AmbKuc06}.
See also \cite{AgrKucQui06} for range results on the spherical mean
transform on 
functions supported in a ball in all dimensions, and \cite{FinRak06} for
range results for the wave trace map for functions supported in the ball in
odd dimensions.

In the work of the first and third authors with Sarah Patch
\cite{FinPatRak04}, several formulas were found to recover a smooth function
$f$ with support in the closure $\overline B$ of the open ball $B \tm \R^n$ 
from the trace of the solution of the wave equation on the product 
$\partial B \times [0,\diam (B)]$ provided that the space dimension is odd.
Specifically, if $u$ is the solution of the initial value problem
\beq{wave}
    u_{\rT\rT} - \Lap u  = 0 ,
    \quad \text{ in } \R^n \times [0,\infty)
\eeq
\beq{novel}
    u(.,\rT=0) = f(.), \quad u_{\rT}(.,\rT=0) =0 ,
\eeq
where $f$ is smooth and has support in the $\overline B,$ then 
several formulas were found to recover $f$ from 
$u(\ps, \rT)$ for $\ps \in S:= \partial B$ and $\rT \in \R^+$.

The first and third authors tried, at that time, to extend the method
to even dimensions, but did not see a way. Recently, the second author
tried numerical experiments using a two dimensional analog of one of
the inversion formulas and found that it gave excellent
reconstructions. This prompted our re-examination of the
problem. Among the results of this paper is a proof of the validity of
this formula.

To describe our results, we introduce some notation.
The spherical mean transform $\M $ is defined by
\begin{equation} \label{spmean}
    (\M f)(x,\rM)
    =
    \frac{1}{\abs{S^{n-1}}}
    \int_{S^{n-1}}
        f(x +\rM\theta)
    \,dS(\theta)
\end{equation}
for $f \in C^\infty(\R^n)$ and $( x, \rM ) \in \R^n \times[0,\infty)$.
In this expression, $\abs{S^{n-1}}$ denotes the area of the unit
sphere $S^{n-1}$ in $\R^n$ and $dS(\theta)$ denotes area measure on
the sphere.  In general, we write the area measure on a sphere of any
radius as $dS$, except when $n=2$ when we write $ds$. 
We will denote the (partial) derivative of a function $q$ with respect
to a variable $\rM$ by $\di_{\rM} q$, except in a few formulas where
the subscript notation $q_{\rM}$ is used. 
At several points we use $D_\rM$ to denote the operator
\[
    (D_{\rM}u) (\rM) :=  \frac{(\di_{\rM} u)(\rM)}{2\rM}
\]
acting on smooth (even) functions $u$ with compact support. 
Moreover, $\ro$ will be used to  denote the   operator  that 
multiplies a function $ u(\rM)$ by $\rM$.

Our first set of results is a pair of inversion formulas for the
spherical mean transform in even dimensions. 
We state and prove these first in dimension two; that is, for the circular
mean transform. 

\begin{theorem} \label{thm:circinv}
    Let $D \subset \R^2$ be the disk of radius $\Rball$ centered
    at the origin, let $S := \di D$ denote the boundary circle, and let
    $f \in C^{\infty}( \R^2 )$ with $\supp f \subset \overline D$.
    Then, for $x \in D$,
    \beq{m-lap}
        f(x) = \frac{1}{2 \pi \Rball} \,
                \Lap_x \int_S \int_0^{2\Rball}
                \rM \, (\M f)(\ps, \rM) \log \wabs{\rM^2 -\abs{ x - \ps}^2 }
                \,d\rM \,ds(\ps)
    \eeq
    and
    \beq{m-inv}
        f(x) = \frac{1}{2\pi \Rball}
               \int_S \int_0^{2\Rball}
               \left( \di_{\rM} \ro \di_{\rM} \M f \right)( \ps, \rM )
                  \log \wabs{ \rM^2 - \abs{x - \ps}^2}
             \,d\rM \, ds(\ps).
    \eeq
\end{theorem}

In Theorem \ref{thm:circinv}, $\di_{\rM} \ro \di_{\rM} \M f $ 
denotes the composition  of $\di_{\rM}$, $\ro$, $\di_{\rM}$ and $\M$ applied
to  $f$.  
The same  convention  will be used throughout the article to denote the 
composition  of any  operators.
 
While $\M f$ has a natural extension to the negative reals as an even
function, we instead take the odd extension in the second variable. 
Then formula \req{m-inv} has the following corollary:

\begin{corollary}\label{2cor}
    With the same hypotheses as in Theorem \ref{thm:circinv}, and 
    $\M f$ extended as an odd function  in the second variable $\rM$,
    $f$ can be recovered for  $x \in D$ by
    \begin{equation}\label{eq:m-inv-hilb}
        f(x) =
            \frac{1}{2\pi \Rball}
                \int_S
                    \int_{-2\Rball}^{2\Rball}
                        \frac{(\ro \di_{\rM} \M f)(\ps,\rM)}{\abs{ x-\ps } -
            \rM } 
                    \, d\rM
                \,ds(\ps),
    \end{equation}
    and
    \begin{equation}\label{eq:filbac}
        f(x) =
        \frac{1}{2\pi \Rball}
        \int_S
            \abs{x-\ps}
                \int_{-2\Rball}^{2\Rball}
                    \frac{(\di_{\rM} \M f)(\ps,\rM)}{\abs{x-\ps} - \rM}
                \, d\rM
        \,
        ds(\ps),
\end{equation}
where the inner integrals are taken in the
principal value sense.
\end{corollary}

These forms are very close to the standard inversion formula for the Radon
transform in the plane \cite[Eq. (2.5)]{Nat86}.

In higher even dimensions we prove a similar pair of results.

\begin{theorem}\label{thm:spherinv}
Let $B \subset \R^n$, $n > 2$ even, be the ball of radius $\Rball$ centered
at the origin, 
let $S := \di B$ be the boundary of the ball, set
\[
    c_n = (-1)^{(n-2)/2}2 ((n-2)/2)! \pi^{n/2} = (-1)^{(n-2)/2}
    [((n-2)/2)!]^2 \abs{S^{n-1}},
\]
and let $f \in C^{\infty}(\R^n)$ have support in $\overline B $.
Then, for $x \in B$,
\begin{align}\label{eq:rhofil_n}
    f(x) 
    &=
    \frac{1}{c_n \Rball} \,
    \Lap_x \int_S \int_0^{2\Rball}
    \log \wabs{\rM^2-\abs{x-\ps}^2}
    (\ro D_{\rM}^{n-2} \ro^{n-2} \M f)(\ps,\rM)
    \,d\rM
    \,dS(\ps),
\\ \label{eq:filback_n}
    f(x) 
    & =
    \frac{2}{c_n \Rball}
    \int_S
        \int_0^{2\Rball}
	    \log\wabs{\rM^2-\abs{x-\ps}^2}
            ( \ro D_{\rM}^{n-1} \ro^{n-1} \di_{\rM} \M f)(\ps,\rM)
        \,d\rM\,
    dS(\ps).
\end{align}
\end{theorem}
Recently, Kunyansky \cite{Kun06} has also established
inversion formulas of the filtered back-projection type for the spherical
mean transform.  His method and results appear to be very different than ours.

For some results, it will be more convenient to use the wave equation
\req{wave} with initial condition \beq{nodis} u(.,\rT=0)= 0,\quad
u_{\rT}(. ,\rT=0) =f(.), \eeq It is obvious that the solution of \req{wave}
with initial values \req{novel} is the time derivative of the solution
of \req{wave} with initial values \req{nodis}.  We denote by $\Po$ the
operator which takes smooth initial data with support in $\overline B$ 
to the solution of \req{wave}, \req{nodis} restricted to $S \times
[0,\infty)$ and by $\W$ the operator taking $f$ to the solution of
\req{wave}, \req{novel} restricted to $S\times [0,\infty)$.  These
operators are simply related by $\W = \di_{\rT} \Po$.
An explicit representation for $\Po$ comes from the well-known formula
\cite{CouHil93} 
\begin{equation}\label{ndsol}
    u(\ps, \rT) =
    \frac{1}{(n-2)!}
    \, \di_{\rT}^{n-2}
    \int_0^\rT
    \rM ( \rT^2 - \rM^2)^{(n-3)/2}
    (\M f)(\ps,\rM) \,d\rM.
\end{equation}
giving the  solution of the initial value problem \req{wave}, \req{nodis},
in dimension  
$n \geq 2$.
We denote by $\Po^\ast$ and $\W^\ast = - \Po^\ast \di_{\rT}$  the formal
$L^2$ adjoints  of $\Po$ and $\W$  mapping  from smooth functions $ u \in
C^\infty(S\times [0, \infty))$  with sufficient decay in the second variable.  
An explicit expression for $\Po^\ast u$   will be given in Section
\ref{sec:wave}. 

We have two types of inversion results for the wave equation. The
first type is based on  the inversion results for the spherical mean
transform, since the spherical mean transform itself can be recovered
from the solution of the wave equation by solving an Abel type
equation. In dimension two, this approach yields the following result.

\begin{theorem}\label{thm:wavefinite}
    Let $D \subset \R^2$  be the open disc with radius $\Rball$
    and let $S := \di D$ denote the boundary circle. 
    Then there exists a kernel function $K: [0,2\Rball]^2 \to \R$ 
    such that for any $f \in C^{\infty}(\R^2)$  with support in $\overline D$
    and any $x\in D$
    \beq{wavefinite}
        f(x)
         =
        \frac{1}{\Rball \pi^2} \,
        \Lap_x
        \int_S \int_0^{2\Rball}
                (\W f)(\ps, \rT) K(\rT, \abs{x-\ps})
        \,d\rT \,ds(\ps).
    \eeq
An analytic expression for $K$ will be given in Section \ref{sec:wave}.
\end{theorem}

Theorem \ref{thm:wavefinite} provides inversion formulas
of the filtered back-projection type for reconstruction of $f$ from  
$(\W f)(\ps, \rT) =( \di_{\rT} \Po f)(\ps, \rT)$ using only data with $\rT \in
[0,2\Rball] $,   
despite the unbounded  support of $\W f$ and $\Po f$ in $\rT.$  

The second type of inversion results holds in all even dimensions and 
takes the following  form.
\begin{theorem}\label{thm:invPW} 
Let $f$ be smooth and supported in closure of the ball $B$ of radius 
$\Rball$ in $\R^{2m}$, and let $\Po f$ and $\W f$ be as above.
Then for $x\in B$
\begin{align} \label{eq:invP}
    f(x)
    &=
        -\frac{2}{\Rball} 
        \left(
        		\Po^\ast \tio \di_{\rT}^2 \Po f
	\right) (x),  
    \\ \label{eq:invW}
    f(x)
    &=
        \phantom{-}
        \frac{2}{\Rball} \left( \W^\ast \tio \W f \right) (x)
        =
        -\frac{2}{\Rball}
        \left(\Po^\ast \di_{\rT} \tio \di_{\rT} \Po f \right) (x).
\end{align}
\end{theorem}

We will prove \req{invP} in dimension $n=2m=2$ directly.  The higher
dimensional case of \req{invP}, and \req{invW} in all dimensions, are
consequences of the following trace identities, relating the $L^2$ inner
product of the initial data to the weighted $L^2$ inner product of the
traces of the solutions of the wave equation.

\begin{theorem}\label{thm:e2m}
    Let $f,g$ be smooth and supported in the ball $B$ of radius $\Rball$,
    in $\R^{2m}$ with $m \geq 1$, let $S:= \partial B$,
    and let $u$ (resp. $v$) be the solution of the initial value problem
    \req{wave},  \req{nodis} with initial value $f$ (resp. $g$).
    Then
    \begin{align}\label{asymm}
        \int_B f(x) g(x) \,dx 
        & =
        -\frac{2}{\Rball}
        \int_S \int_0^{\infty} 
        		\rT u_{\rT\rT} (\ps, \rT) v(\ps, \rT)
        \,d\rT \,dS(\ps),
    \\  \label{symm}
    \int_B f(x) g(x)  \,dx 
    &=
    \phantom{-} 
    \frac{2}{\Rball} 
    \int_S \int_0^{\infty} 
    	\rT u_{\rT}(\ps, \rT) v_\rT(\ps, \rT)
    \,d\rT\,dS(\ps).
\end{align}
\end{theorem}

In the proof of this theorem, (\ref{asymm}) for $n=2$ follows from
\req{invP} for $n=2,$ while (\ref{asymm}) in higher even dimensions is derived
from the $n=2$ case; (\ref{symm}) is a consequence of (\ref{asymm}) in all
dimensions.
We remark that these identities were already proved in \cite{FinPatRak04} for
odd dimensions, and so they hold for all dimensions.

Section 2 is devoted to the proof of the inversion formulas for the
spherical mean transform, that is, Theorems 1, 2, and Corollary 1. Section 3
treats the wave equation and contains the proofs of Theorems 3, 4, and
5. This is 
followed by a section reporting on the implementation of the various
reconstruction formulas of the preceding sections and results of
numerical tests, in dimension two.

\section{Spherical Means}

In this section we prove the Theorems related to the inversion from
spherical means 
and Corollary \ref{2cor}. We begin by establishing an elementary integral
identity, which is the key to the results in this paper.

\begin{pr} \label{pr:keykeyident}
Let $D \tm \R^2$ be the disk of radius $\Rball$, and let $S=\di D$ be the
boundary circle. Then for $x$, $y \in D$ with $x \neq y$,
\beq{keyident}
    \int_S \log \left| \vert x-\ps\vert^2 -\vert y -\ps\vert^2 \right| \,
    ds(\ps) 
    =
    2\pi \Rball \log\vert x-y\vert +2\pi \Rball\log\Rball.
\eeq
\end{pr}

\begin{proof}
Let $x\neq y$ both lie in $ D$ and let $I$ denote the integral
on the left on \req{keyident}. Expanding the argument of the
logarithm  as
\beqn
    \wabs{\abs{x - \ps }^2 - \abs{ y - \ps }^2} =
    2 \Rball \abs{x - y} \wabs{\left(
    \frac{ x + y }{2\Rball}-\frac{\ps}{\Rball}\right) \cdot
    \frac{ x - y }{ \abs{ x - y}}},
\eeqn
setting $e:= \frac{x-y}{\abs{x-y}}$, and writing $\ps=\Rball\theta$ for 
$\theta\in S^1$, we have
\begin{equation}\label{eq:help0}
    I =
    2\pi \Rball \log \left( 2\Rball\abs{x-y}  \right)
    + \Rball \int_{S^1}
        \log \abs{e\cdot \theta -a} \,d\theta,
\end{equation}
where \[ a = \frac{x+y}{2\Rball}\cdot e =
\frac{\abs{x}^2-\abs{y}^2}{2\Rball\abs{x-y}}.\]
We note that $\abs{a} < 1.$

Using the parameterization $\theta = \cos(\phi) e + \sin(\phi) e^\bot,$
the integral term on the right of \req{help0} 
has the form
\[ \Rball \int_0^{2\pi} \log\abs{\cos\phi -a}\,d\phi.\]
Writing $a = \cos \alpha$ and using the \textit{sum to product}
trigonometric identity 
$\cos \phi -  \cos \alpha  = -2\sin \left( (\phi +\alpha)/2 \right)\sin
\left( (\phi -\alpha)/2 \right)$, 
this is equal to
\[
    \Rball
    \int_0^{2\pi} \left(
    \log 2 + \log \abs{ \sin \left( (\phi +\alpha)/2 \right) } +
    \log \abs{ \sin \left( (\phi-\alpha)/2 \right)}
     \right)  d\phi.
\]
By periodicity, and two linear changes of variable, this reduces to
\[    \Rball  \int_0^{2\pi} \left( \log 2 +2\log\vert \sin(\phi/2)\vert \right)
    d\phi =
    2 \Rball \pi \log 2
    + 4\Rball \int_0^{\pi} \log\sin u \,du, \]
which is independent of $\alpha$, and hence of  $x$ and $y$.
The latter integral in can be found in tables,
and is equal to $-\Rball \pi \log2$,  so the sum is $-2\pi\Rball \log 2$.
Substituting in \req{help0}
gives the desired result.
\end{proof}

Proposition \ref{pr:keykeyident} is already enough to establish
Theorem \ref{thm:circinv}.

\begin{proof1}
Let $f\in C^\infty(\R^2)$ be supported  in $\overline D$ and let $\ps$ be any
point in $S = \di B$.
Using the definition of  $\M f$ and  Fubini's theorem, we have that
\beq{pullback}
    \int_0^{2\Rball}
    ( \ro \M f)(\ps,\rM) q(\rM)\,d\rM
    =
    \frac{1}{2\pi}\int_{\R^2}
    f(\ps +z) q( \abs{z} ) \,dz,
\eeq
for any measurable function $q$ provided that the product of functions
on the right is absolutely integrable.
Applying this with $q(\rM) = \log \wabs{ \rM^2 -\abs{x-\ps}^2 }$ and making the
change of variables $y = \ps + z$ gives
\begin{multline*}
    \int_S
        \int_0^{2\Rball}
             (\ro \M) (f)(\ps,\rM)
            \log \wabs{ \rM^2 -\abs{x-\ps}^2 }
            \,d\rM\, ds(\ps) \\
     = \frac{1}{2\pi} \int_S
    \int_{\R^2} f(y)
    \log \wabs{ \abs{y-\ps}^2 -\abs{ x-\ps}^2}
    \,dy \, ds(\ps).
\end{multline*}
Fubini's theorem again justifies the change of order of integration in
the iterated integral on the right hand side, and so 
\[
    \frac{1}{2\pi}
    \int_{\R^2}  f(y) \int_S
    \log \wabs{ \abs{ y-\ps}^2  - \abs{ x-\ps}^2 } \,ds(\ps) \,dy
        =
    \frac{2 \pi \Rball }{2\pi}
    \int_{\R^2}
        f(y) ( \log \abs{x-y} + \log \Rball )
    \,dy 
\]
upon application of \req{keyident}.  Recalling that for any constant $c$,
$1/(2\pi) \log \abs{x-y}+c $ is a fundamental solution of the Laplacian in
$\R^2$, we have
\[
    f(x) =
    \frac{1}{2\pi \Rball} 
    \, \Lap_x
    \int_S
        \int_0^{2\Rball}
             (\ro \M f)(\ps,\rM) \log \abs{\rM^2 -\abs{x-\ps}^2}
    \,d\rM\,ds(\ps),
\]
which proves \req{m-lap}.

The second formula, \req{m-inv}, has a similar proof.  In this case,
we use that the spherical means satisfy the Euler-Poisson-Darboux equation
\cite{CouHil93}
\[
    (\di_{\rM}^2 \M f)(x,\rM) +\frac{1}{\rM}
    (\di_{\rM} \M f)(x,\rM) = (\Lap \M f)(x,\rM) = (\M \Lap f)(x,\rM).
\]
The left hand side of the  Darboux equation may be written as   $(1/\rM)(
\di_\rM \rM \di_\rM \M f)(x, \rM) $, so 
the expression on the right of \req{m-inv} may be
rewritten as 
\beq{help1}
   \frac{1}{2\pi \Rball}\int_S
    \int_0^{2\Rball} (\ro \M \Lap f)(\ps,\rM)
    \log \wabs{ \rM^2 -\abs{x-\ps}^2 }
    \,d\rM\,ds(\ps).
\eeq
Again applying \req{pullback}, now with the function $q(\rM) = \rM \log
\wabs{ \rM^2 -\abs{x-\ps}^2 }$ 
and $\Lap f $ instead of $f$, interchanging the order of integration and
using \req{keyident}  
shows that the expression \req{help1} is equal to  
\beqn
   \frac{1}{2\pi}  \int_{\R^2}
            \Lap_y
                f(y)
                (
                \log \vert x - y \vert
                + \log \Rball
                )\,dy
    =
    f(x),
\eeqn
since no boundary terms arise in view of the support hypothesis on $f$. \eproof
\end{proof1}

\begin{proof2}
Let $x \in D$,  and let
\beqn
U(\ps, x)
:=
\int_0^{2\Rball}
     \left( \di_{\rM}  \ro \di_{\rM} \M f\right)( \ps, \rM )
    \log \wabs{ \rM^2 - \abs{x - \ps}^2}
\,d\rM
\eeqn
denote the inner integral in \req{m-inv}.
Taking the  support of $f$ into  account, writing the logarithm as
\beqn
    \log \wabs{ \rM^2 -\abs{x-\ps}^2 } =
    \log \wabs{ \rM - \abs{ x-\ps} } +
    \log \wabs{ \rM +\abs{x-\ps}},
\eeqn
and integrating \req{m-inv} by parts leads to
\beqn
  U(\ps, x)
    &= -P.V. \int_0^{\infty} \frac{(\ro \di_{\rM} \M f)(\ps, \rM)}{\rM -
    \abs{x-\ps}}\,d\rM - \int_0^{\infty} \frac{(\ro \di_{\rM} \M f)(\ps,
    \rM)}{\rM + \abs{x-\ps}}\,d\rM.
\eeqn
Here we have used that the distributional derivative of $\log \abs{r}$ is
$P.V. \,\frac{1}{r}$ as well as an ordinary integration by parts.  
Therefore \req{m-inv} implies
\beq{cor-help}
    f(x)
    &
    = \frac{1}{2\pi \Rball} \int_S U(\ps, \rM) ds(\ps)
    \\
    &
    =
    \frac{-1}{2\pi \Rball}
    \int_S \int_0^{2\Rball}
        \frac{(\ro \di_{\rM} \M f)(\ps,\rM)}{\rM -\abs{x-\ps}}
     \,d\rM\,ds(\ps)
   + \frac{-1}{2\pi \Rball}
    \int_S \int_0^{2\Rball}
       \frac{(\ro \di_{\rM} \M f)(\ps,\rM)}{\rM +\abs{x-\ps}}
    \,d\rM\,ds(\ps),
\eeq
where the inner integral of the first term on the right is taken in the
principal value sense. 
The odd extension of $\M f$, $\M f (\ps, -\rM) := -\M f (\ps, \rM) $, is
smooth on $\R$ since $\M f$ vanishes to infinite order at $\rM = 0$ by the
support hypothesis on $f$ 
and  $(\ro \di_{\rM} \M f)(\ps,\rM)$ is an odd function in $\rM.$
Substituting $\rM = -\rM$ in the second  integral in \req{cor-help} gives
\beqn
    f(x)
    =
    \frac{-1}{2\pi \Rball}
        \int_S \int_{0}^{2\Rball}
            \frac{(\ro \di_{\rM} \M f)(\ps,\rM)}{\rM - \abs{ x-\ps }}
        \,d\rM\,ds(\ps)
    + \frac{-1}{2\pi \Rball}
    \int_S \int_{-2\Rball}^{0}
      \frac{(\ro \di_{\rM} \M f)(\ps,\rM)}{\rM -\abs{x-\ps}}
      \,d\rM\,ds(\ps)
\eeqn
and hence
\[
    f(x) =
    \frac{1}{2\pi \Rball}
    \int_S \int_{-2\Rball}^{2\Rball}
        \frac{(\rM \di_{\rM} \M f)(\ps,\rM)}{\abs{x-\ps} - \rM }
    \,d\rM
    \,ds(\ps).
\]
This is (\ref{eq:m-inv-hilb}). To prove (\ref{eq:filbac}), it suffices to write
\[
    \frac{\rM}{\abs{x-\ps} - \rM}
    =
    -1 + \frac{\abs{x-\ps}}{\abs{x-\ps} - \rM}
\]
in (\ref{eq:m-inv-hilb}) and to note that
$\int_{-2\Rball}^{2\Rball} (\di_{\rM} \M f) (\ps,\rM) \,d\rM =0$,
by the support hypothesis on $f$. \eproof

\end{proof2}

\subsection{Proof of Theorem \ref{thm:spherinv}}

We have found several proofs of Theorem \ref{thm:spherinv}, the extension of
Theorem \ref{thm:circinv} to higher even dimensions. The one we present is
based on reduction of the higher dimensional problem to the two dimensional
case already established. Another, which is not presented in this
article, is based on an extension of \req{keyident} to higher dimensions.

We first observe that by a dilation, we may reduce the problem to the case
when $f$ is supported in the unit ball.  Tracing through the formulas
\req{rhofil_n} and \req{filback_n} it is routine to verify that scaling from
the unit ball to the ball of radius $\Rball$ introduces a factor of
$\Rball$.  To simplify notation, we shall now suppose that $f$ is supported
in the unit ball $B$. Let $Q$ and $N$ denote the operators
 \begin{align} \label{eq:Q}
    (Q f) (x) 
    &=
    \Lap_x \int_S \int_0^{2}
    (\ro D_{\rM}^{n-2} \ro^{n-2} \M f)(\ps,\rM)
    \log \wabs{\rM^2-\abs{x-\ps}^2}
    \,d\rM
    \,dS(\ps),
\\ \label{eq:N}
    (N f)(x) 
    & =
    \int_S
        \int_0^{2}
            ( \ro D_{\rM}^{n-1} \ro^{n-1} \di_{\rM} \M f)(\ps,\rM)
            \log\wabs{\rM^2-\abs{x-\ps}^2}
        \,d\rM\,
    dS(\ps) ,
\end{align}
that map $f \in C^\infty(\R^{n}) $ supported in $\overline B$ to constant
multiples of the  
the right hand sides  of \req{rhofil_n} and \req{filback_n}. Moreover 
$\langle f, g\rangle$ denotes the $L^2$ product of two functions supported in 
$\overline B$.  To establish $Q f = c_n f$ and $N f = (c_n/2) f$ we will use  
the following auxiliary results.
 
\begin{pr}\label{pr:NQ}
Let $f$, $g$ be smooth and supported in $\overline B$. 
Then 
\begin{equation} \label{eq:NQ} 
	\int_{\R^n} (Qf) (x) g(x) \, dx
	= 
	\langle Qf, g \rangle 
	= 
	2 \langle f, Ng \rangle
	=
	2 \int_{\R^n} f (x) (N g) (x) \, dx.
\end{equation}
\end{pr}

\begin{proof}
Let $F= \M f$  and $G = \M g$.
Using the self-adjointness of $\Lap$, applying Fubini's theorem and  an 
$n$-dimensional analogue 
of \req{pullback},  we obtain
\begin{align}
    \langle Qf, g \rangle
    &=
    \int_B 
    	\left(
		\int_S 
		\int_0^2
			(\ro D_{\rM}^{n-2} \ro^{n-2} F)(\ps,\rM)
    			\log \wabs{\rM^2-\abs{x-\ps}^2}
    			\, d\rM
		\,dS(\ps)
	\right)
    (\Lap_x g)(x) \,dx \notag
    \\
    &=
    \abs{S^{n-1}}
    \int_S 
	\int_0^2 
		\left(
			\int_0^2
			(\ro D_{\rM}^{n-2} \ro^{n-2} F)(\ps,\rM)
    			\log \wabs{\rM^2-\rH^2}
    			(\M \Lap_x g) (\ps, \rH) \rH^{n-1}
		\, d\rH
		\right)
	\,d \rM 
    \,dS(\ps) \notag
\\
    &=
    \abs{S^{n-1}}
    \int_S 
	\int_0^2 
		\left(
			\int_0^2
			(\ro D_{\rM}^{n-2} \ro^{n-2} F)(\ps,\rM)
    			\log \wabs{\rM^2-\rH^2}
		\, d\rM
		\right)
	(\M \Lap_x g) (\ps, \rH) \rH^{n-1}
	\,d \rH
    \,dS(\ps)\notag 
\\
   &= \abs{S^{n-1}}\int_S \int_0^2\left(
           \int_0^2(\ro D_{\rM}^{n-2} \ro^{n-2} F)(\ps,\rM)
    			\log \wabs{\rM^2-\rH^2}
		\, d\rM
		\right)
         \di_{\rH}\rH^{n-1}\di_{\rH}G(\ps,\rH)\,d\rH\,dS(\ps).\label{eq:aux2}  
\end{align}
To justify the last equation it is used that $G$ satisfies the
Euler-Poisson-Darboux equation and the identity $\rH^{n-1}(\di_{\rH}^2
+ \frac{n-1}{\rH} \di_{\rH}) = \di_{\rH}(\rH^{n-1}\di_{\rH}).$
Applying the identities $(D_\rM^{n-2})^* \rM \log\abs{\rM^2 - \rH^2} =
(-1)^{n-2} \rM D_\rM^{n-2} \log\abs{\rM^2-\rH^2} = \rM D_{\rH}^{n-2}
\log\abs{\rM^2-\rH^2}$ in two stages to the last expression, this
becomes
\begin{align*}
    \abs{S^{n-1}}
    &\int_S
    	\int_0^2
    		\left(
			\int_0^2 
			\rM^{n-1} F(\ps,\rM)
    			D_{\rH}^{n-2} \log\abs{\rM^2-\rH^2}\, d\rM
		\right)
    	\left( \di_{\rH} \rH^{n-1}\di_{\rH} G(\ps,\rH)\right)
    \,d\rH \,dS(\ps) 
  \\
  &=
    \int_S
    	\int_0^2
    		\left(
			\int_B 
			f(y)
    			D_{\rH}^{n-2} \log\abs{\abs{y-\ps}^2-\rH^2}\, dy
		\right)
    	\left(\di_{\rH} \rH^{n-1} \di_{\rH}G(\ps,\rH)\right) \,d\rH \,dS(\ps)
\\
     &=	\int_B 
	\left( 
		\int_S \int_0^2  
		\log \wabs{\abs{y-\ps}^2 - \rH^2}
		( (D_{\rH}^{\ast})^{n-2} \di_{\rH} \rH^{n-1} \di_{\rH}
	G)(\ps,\rH)  
	\,d\rH\,dS(\ps)
	\right)
	f(y) \,dy
\end{align*}   
after applying Fubini's theorem. This is finally seen to be equal to 
$ \langle f, 2 N g\rangle$  since $(D_{\rH}^{n-2})^\ast \di_{\rH} = 2\rH
(-1)^{n-2} D_{\rH}^{n-1}$. 
\end{proof}

We now look at  the spherical means of  products
\beq{fPol}
    f(x) &= \rP^k  \alpha(\rP) \Phi(\theta),
\eeq
where  $x = \rP \theta$ with $\rP \geq 0$, $\theta \in S^{n-1}$, 
$\Phi$ is a spherical harmonic of degree $k$, and $\alpha: \R \to \R$ is an
even smooth function supported in $[-1,1]$. 
Let $F := \M f$ be extended to an even function in the second component and
let  $\nu = n + 2k$.   
Then $F$ satisfies the initial value problem (IVP) for the
Euler-Poisson-Darboux equation
\begin{align} \label{eq:Fpde}
    \left(\di^2_{\rM}F  + \frac{n-1}{\rM} \di_{\rM}F \right)(x, \rM)
    =
    \Lap_x F (x, \rM),
    \qquad 
    & (x, \rM) \in \R^{n} \times \R
\\
    F( x, 0)
    =
    \alpha(\rP) \rP^k \Phi(\theta),
    \quad
    \di_{\rM}F(x, 0) = 0,
    \qquad 
    & x \in \R^n,
\label{eq:Fic}
\end{align}
and, conversely, any solution of \req{Fpde}, \req{Fic}, is the spherical
mean of the initial values.
The unique solution of \req{Fpde}, \req{Fic} has the form
$F(x  , \rM) = \rP^{k}  A(\rP, \rM) \Phi(\theta)$  where $A(\rP,\rM)$ is the
solution of 
the  IVP
\begin{align}
     (L_n A)(\rP,\rM)
    = \left(\di^2_{\rP}A + \frac{\nu-1}{\rP} \di_{\rP}A \right)(\rP,\rM),
    \qquad 
    & (\rP,r) \in \R^2,
\label{eq:Apde}
\\
    A(\rP, 0)=\alpha(\rP),
    \quad
    \di_{\rP}A(\rP, 0) = 0,
    \qquad
    & \rP \in \R.
\label{eq:Aic}
\end{align}
Here  $ (L_n A)  (\rP, \rM) :=  ( \di_{\rM}^2 A + \frac{n-1}{\rM} \di_{\rM}A)
(\rP,\rM)$. 

We recall that  the operator $D_{\rM}$ satisfies $ L_n D_{\rM} = D_{\rM}
L_{n-2}$ and 
for any $\mu \in \N$
\[
    \left(
    \di_\rM^2 +
    \frac{1-\mu}{\rM} \di_{\rM}
    \right)
    ( \rM^\mu w)  =
    \rM^\mu
    \left(
        \di_\rM^2 +  \frac{1+\mu}{\rM} \di_{\rM}
    \right) w,
\]
that is $L_{2-\mu} \rM^\mu = \rM^{\mu} L_{\mu+2}$.  So
\beq{Ldescent}
    ( L_{2-\mu+ 2 \sigma} D_{\rM}^\sigma \rM^\mu w) (\rM)
    =
    (D_{\rM}^{\sigma} L_{2-\mu} \rM^\mu w ) (\rM)
    =
    ( D_{\rM}^{\sigma} \rM^\mu L_{\mu+2} w ) (\rM).
\eeq
If we set  $\mu = n-2$ and $\sigma = (n - 2)/2$ in \req{Ldescent},
then $\mu+2 = n$ and $2 - \mu + 2 \sigma = 2$.  Therefore
\beq{LDg}
    ( L_{2}  D_{\rM}^{(n-2)/2} \rM^{n-2} w) (\rM)
    =
    ( D_{\rM}^{(n-2)/2} \rM^{n-2} L_n w) (\rM).
\eeq

Now we set  
\beq{Hdef}
H(\rP,\rM) := \dfrac{1}{ (\ntw)!} (D_{\rM}^{(n-2)/2} \rM^{n-2}
A) (\rP,\rM).
\eeq 
Since $A(\rP,\rM)$ is even in $\rM$ and $D_{\rM}$ corresponds to
differentiation with   
respect to $\rM^2$, $H(\rP,\rM)$ is even in $\rM$. 
Moreover, by \req{Aic}, 
$H(\rP,0) = \frac{1}{(\ntw)!} A(\rP,0) ( D_{\rM}^{(n-2)/2}
\rM^{n-2} ) = 
\alpha(\rP)$, 
and therefore from \req{Apde} and \req{LDg}  it follows that $H$ 
is the solution of  the IVP
\begin{align}
    \left(
    \di^2_{\rM}H  + \frac{1}{\rM} \di_{\rM}H
    \right)
    (\rP,\rM)
    =
    \left(\di^2_{\rP}H  + \frac{\nu-1}{\rP}\di_{\rP}H \right)(\rP,\rM),
    \qquad 
    & (\rP,\rM) \in \R^2,
\label{eq:Hpde}\\
    H(\rP,0)= \alpha(\rP),  \quad
    \di_{\rM}H(\rP,0) = 0,  
    \qquad 
    & \rP \in \R,
\label{eq:Hic}
\end{align}

\begin{pr}\label{pr:intID}
Let $A_i(\rP,\rM)$, $i=1,2$ solve \req{Apde} with $n=2,$ subject to initial 
conditions $A_i(\rP,0) = \alpha_i(\rP),$
$\di_{\rM}A_i(\rP,0) =0,$ where $\alpha_i$ are smooth even functions with
support in $[-1,1]$ and $\nu \geq 2$ is even. Then
\beq{IDnu}
\int_0^1 \rP^{\nu -1} \alpha_1(\rP) \alpha_2(\rP)\,d\rP = -\int_0^2\int_0^2
\rM A_1(1,\rM) \di_{\rH}\log\abs{\rM^2-\rH^2}\, \rH (\di_{\rH}A_2)(1,\rH)
\,d\rH\,d\rM. 
\eeq
\end{pr}
\begin{proof}
Let $k = (\nu -2)/2$ and let $\Phi(\theta)$ be a nontrivial real
circular harmonic of degree $k$. Then $F_i(x,\rM) := A_i(\rP,\rM)
\rP^k \Phi(\theta)$ satisfies \req{Fpde}, \req{Fic} for $n=2$, and
so is the circular mean of its initial value, $f_i(x) = \alpha_i(\rP)
\rP^k \Phi(\theta).$ By \req{m-lap}, $f_1 = \frac{1}{2\pi} Q f_1,$ and using
\req{aux2} gives
\beqn 
\langle f_1, f_2\rangle &= \frac{1}{2\pi} \langle Q f_1,
f_2\rangle\\ &= \int_S \int_0^2\int_0^2 r F_1(\ps,\rM) \log\abs{\rM^2
-\rH^2} (\di_{\rH} \rH\di_{\rH}
F_2)(\ps,\rH) \,d\rM\,d\rH\,ds(\ps).  
\eeqn
Taking account the form of $F_i$ and that $\rP =1$ on $S$, this may be 
rewritten as 
\beq{cyl}
\langle f_1, f_2\rangle = \int_S \Phi^2(\ps)\, ds(\ps) \int_0^2 \int_0^2 \rM
A_1(1,\rM) \log\abs{\rM^2 -\rH^2} (\di_{\rH} \rH \di_{\rH} 
A_2)(1,\rH) \,d\rM \,d\rH.
\eeq
Appealing to the form of $f_i = F(x,0)$, 
\beq{iniV}
\langle f_1, f_2\rangle &= \int_0^1 \rP (\rP^k \alpha_1) (\rP^k \alpha_2) 
\,d\rP \int_S \Phi^2(\ps) \,ds(\ps)\\
&= \int_0^1 \rP^{\nu -1} \alpha_1(\rP) \alpha_2(\rP) \,d\rP 
\int_S \Phi^2(\ps) \, ds(\ps).
\eeq
Since $\int_S \Phi^2(\ps)\, ds(\ps) \neq 0$, a comparison of  \req{cyl} and 
\req{iniV} and an integration by parts on the right side of \req{cyl}
establishes \req{IDnu} which completes the proof. 
\end{proof}

\begin{proof3}
Let $\{\Phi_j\}$ be an orthonormal basis for the spherical harmonics on
$S^{n-1}$, and consider $f_i,$ $i=1,2$ of the form \req{fPol} with $\alpha =
\alpha_i$ and $\Phi = \Phi_{j_i}$
of possibly different degrees. Let $F_i $ be the even extensions of $\M f_i$
as above. Then by orthogonality, $\langle f_1, f_2\rangle =0$ unless $j_1 =
j_2$, in which case 
\beq{f12ini}
\langle f_1, f_2\rangle = \int_0^1 \rP^{\nu -1} \alpha_1(\rP)\alpha_2(\rP)
\,d\rP,
\eeq
with $\nu = n +2k,$ where $k$ is the degree of $\Phi_{j_1}$.
Evaluating $\langle Qf_1, f_2\rangle$ by \req{aux2}, and using that $F_i =
\rP^{k_i} A_i(\rP,\rM) \Phi_{j_i},$ we see that it is also
zero unless $j_1 = j_2.$ In this case we have
\beq{sdegA}
\langle Q f_1&, f_2\rangle = \abs{S^{n-1}}\int_0^2 \int_0^2 (\ro
D_{\rM}^{n-2} \ro^{n-2} A_1)(1,\rM) \log\abs{\rM^2 -\rH^2}
(\di_{\rH}\rH^{n-1}\di_{\rH}A_2)(1,\rH)\, d\rM\,d\rH\\
&= \abs{S^{n-1}}\int_0^2 \int_0^2
(D_{\rM}^{\frac{n-2}{2}}\ro^{n-2}A_1)(1,\rM)(D_{\rM}^*)^{\frac{n-2}{2}}
(r\log\abs{\rM^2-\rH^2}) (\di_{\rH}\rH^{n-1}\di_{\rH}A_2)(1,\rH)
\,d\rM\,d\rH\\
&= \abs{S^{n-1}}\int_0^2 \int_0^2 \ro (D_{\rM}^{\frac{n-2}{2}} \ro^{n-2} A_1)(1,\rM)
D_{\rH}^{\frac{n-2}{2}}
\log\abs{\rM^2-\rH^2}(\di_{\rH}\rH^{n-1}\di_{\rH}A_2)(1,\rH) \,d\rM\, d\rH,
\eeq
since $(D_\rM^{(n-2)/2})^* \rM \log\abs{\rM^2 - \rH^2} =
(-1)^{(n-2)/2} \rM D_\rM^{n-2} \log\abs{\rM^2-\rH^2} = \ro D_{\rH}^{(n-2)/2}
\log\abs{\rM^2-\rH^2}$. Applying the adjoint (distributional derivative)
again in \req{sdegA},
\beqn
\langle Q f_1&, f_2\rangle = \abs{S^{n-1}} \int_0^2 \int_0^2 \ro
(D_{\rM}^{\frac{n-2}{2}} \ro^{n-2} A_1)(1, \rM) \log\abs{\rM^2 -\rH^2}
 (D_{\rH}^{\frac{n-2}{2}})^*(\di_{\rH}\rH^{n-1} \di_{\rH}A_2)(1,\rH)
\, d\rM \, d\rH \\
&= \abs{S^{n-1}}\int_0^2 \int_0^2 \ro (D_{\rM}^{\frac{n-2}{2}}
\ro^{n-2}A_1)(1,\rM) \log\abs{\rM^2 -\rH^2} (-1)^{\frac{n-2}{2}} (\di_{\rH}
D_{\rH}^{\frac{n-2}{2}} \rH^{n-1} \di_{\rH}A_2)(1,\rH)\, d\rM \, d\rH\\
&= \abs{S^{n-1}}(-1)^{n/2} \int_0^2 \int_0^2 \ro (D_{\rM}^{\frac{n-2}{2}}
\ro^{n-2} A_1)(1,\rM) \di_{\rH}\log\abs{\rM^2-\rH^2} (D_{\rH}^{\frac{n-2}{2}}
\rH^{n-1} \di_{\rH} A_2)(1,\rH)\, d\rM \, d\rH.
\eeqn
We now use the following identity, which is readily proved by induction,
\[ D_{\rH}^{(n-2)/2} \rH^{n-1} \di_{\rH} q = \rH \di_{\rH} D^{(n-2)/2}_{\rH}
\rH^{n-2} q\] taking $q = A_2,$ and observe that defining $H_i$ by
\req{Hdef} with $A = A_i$, it then holds that 
\beqn
\langle Q f_1, f_2\rangle &= \abs{S^{n-1}}(-1)^{n/2} [((n-2)/2)!]^2 \int_0^2
\int_0^2 \ro H_1(1,\rM) \di_{\rH}\log\abs{\rM^2-\rH^2} \,\rH (\di_{\rH}
H_2)(1,\rH) \,d\rM \,d\rH. 
\eeqn
By \req{Hpde} and \req{Hic} the $H_i$ satisfy the hypotheses of Proposition
\ref{pr:intID}  with initial data $\alpha_i$, so by \req{IDnu} the 
expression on the right is equal to 
\[ (-1)^{(n-2)/2}\abs{S^{n-1}}[((n-2)/2)!]^2 \int_0^1 \rP^{\nu -1}
\alpha_1(\rP) \alpha_2(\rP) \, d\rP.\]
Thus we have proved that for $f_i$ of the form above
\beq{iniID}
\langle Q f_1 ,f_2\rangle = (-1)^{(n-2)/2}\abs{S^{n-1}} [((n-2)/2)!]^2
\langle f_1, f_2\rangle.
\eeq
We note that the constant on the right  is $c_n$ of Theorem
 \ref{thm:spherinv}.
By linearity and orthogonality of spherical harmonics, this still holds when
either $f_1$ or $f_2$ is replaced by a finite linear combination of such
functions. The set of finite  linear combinations of functions of form
\req{fPol} is dense in $L^2,$ and so we have $Q f  = c_n f$ in $L^2$ when
$f$ is a finite linear combination of functions of the form \req{fPol}.
Now, let $g$ be smooth with with support in the unit ball.  Applying
Proposition \ref{pr:NQ}, it 
follows that
\begin{equation}\label{eq:adj}
    \langle f, N g\rangle =(1/2) \langle Q f, g\rangle = (c_n/2) \langle f,
    g\rangle.
\end{equation}
for all $f$ as above.   Since \req{adj} holds for a
dense subset of functions $f$ in $L^2(B)$, it implies that $N g = (c_n/2) g$
almost everywhere in $B$.  However, $N g$ is easily seen to be a continuous
function, and so $N g = (c_n/2) g$ holds pointwise in $B$, which is
\req{filback_n}. But if $N$ is a multiple of the identity, then so $Q$, and
the proof is complete.
\eproof
\end{proof3}

\section{The Wave Equation}
\label{sec:wave}
We begin the analysis of recovery of initial data from the trace of the
solution of the wave equation on the lateral boundary of the cylinder.
As mentioned in the introduction, we have two types of inversion
results. The first, Theorem \ref{thm:wavefinite}, is really a corollary of
one of the 
inversion formulas for circular means from the previous section.

\begin{proof4}
Let $u(x,\rT)$ to be  the solution of the IVP \req{wave}, \req{novel} in
dimension two. Then by (\ref{ndsol}),   
\[
	u(\ps,\rT) = \di_{\rT} \int_0^t \frac{( \ro \M f)(\ps,\rM)
	}{\sqrt{\rT^2 -\rM^2}} \, d\rM. 
\]
for $\ps \in S$. We can recover the circular means from $u$ by the standard
method of 
inverting an Abel type equation. The details are not hard and may be found,
for example, in \cite{Nat86}. 
The result is
\beq{p2m}
	 (\M f)(\ps,\rT)
	=
	\frac{2}{\pi}
	\int_0^{\rM} \frac{  u(\ps,\rT)}{\sqrt{\rM^2 - \rT^2}}  \,d\rT .
\eeq
Inserting \req{p2m} into the inversion formula \req{m-lap} for $\M$ and 
applying Fubini's theorem, gives, for $x \in D,$
\begin{align*}
f(x) 
& = 
\frac{1}{2\pi \Rball} 
\, \Lap
\int_S 
	\int_0^{2\Rball} 
		(\ro \M f)(\ps,\rM)
		\log \wabs{\rM^2 -\abs{x-\ps}^2} 
	\,d\rM
\,ds(\ps).
\\
& = 
\frac{1}{\Rball \pi^2} 
\, \Lap 
\int_S 
	\int_0^{2\Rball} 
		\rM 
		\int_0^{\rM}
			\frac{u(\ps,\rT)}{\sqrt{\rM^2 - \rT^2}}  
			\log \wabs{ \rM^2 -\abs{x-\ps}^2 }
		\,d\rT
	\,d\rM
\,ds(\ps) 
\\
&= 
\frac{1}{\Rball \pi^2} 
\, \Lap
\int_S 
	\int_0^{2\Rball} 
		u(\ps,\rT) 
		\int_{\rT}^{2\Rball}
			\frac{\rM}{\sqrt{\rM^2-\rT^2}}
			\log\wabs{\rM^2 -\abs{x-\ps}^2} 
		\,d\rM
	\,d\rT
\,ds(\ps)\\
	&=   
\frac{1}{\Rball\pi^2}
  \,   \Lap	\int_S 
		\int_0^{2\Rball} 
			u(\ps,\rT)
			K(\rT,\abs{x-\ps}) 
		\,d\rT
	\,ds(\ps).
\end{align*}
Since  $u(\ps,\rT) = (\W f)(\ps,\rT),$ this is \req{wavefinite},
with
\beq{K}
K(\rT, \rH)
:=
\int_{\rT}^{2\Rball}
\frac{\rM}{\sqrt{\rM^2-\rT^2}}
\log\wabs{\rM^2 - \rH^2} d\rM.
\eeq 
The integral in \req{K} can be evaluated exactly. 
 For the sake of
completeness,  we give the analytic expression.  If we
substitute  $\rM = \sqrt{\rT^2 + \xi^2}$ in \req{K}, then $ d \rM =
(\xi/\rM) d\xi$ and  
thus
\beqn
K(\rT, \rH) 
&= 
\int_{0}^{\sqrt{4\Rball^2-\rT^2}}
	\log \wabs{ \xi^2 + (\rT^2 - \rH^2) } d\xi
\\&
= 
\sqrt{4\Rball^2-\rT^2}
\left( 	- 2 + \log \abs{4 \Rball^2 - \rH^2 }
\right) +\Gamma(\rT,\rH),
\eeqn
where
\beqn \Gamma(\rT,\rH) = \begin{cases}\sqrt{\rH^2
      -\rT^2}\log\frac{\sqrt{4\Rball^2-\rT^2} + \sqrt{\rH^2
	-\rT^2}}{\sqrt{4\Rball^2 -\rT^2}-\sqrt{\rH^2-\rT^2}} & 
      \rT < \rH, \\2 \sqrt{\rT^2 - \rH^2} 
\arctan \sqrt{ \frac{4\Rball^2-\rT^2}{\rT^2 - \rH^2}}  & \rT >
  \rH.\end{cases} 
\eeqn
\eproof
\end{proof4}

For the second type of inversion formula, we start by
deriving a representation of the formal adjoint $\Po^\ast$, for $n=2$.
For any continuous function $G(\ps,\rT)$ on $S \times [0,\infty)$ that has
a small amount of decay as $\rT \to \infty$, by Fubini's theorem,  we have
\begin{align*}
    \langle \Po f,\, G \rangle
    & =
    \int_S
    \int_0^\infty (\Po f)(\ps,\rT)
     \, G(\ps,\rT) \, d\rT \, ds(\ps) \\
    & =
    \frac{1}{2 \pi} \int_S \int_0^\infty
    \, G(\ps,\rT)
    \left( \int_0^\rT \frac{\rM}{\sqrt{\rT^2 - \rM^2}}
    \int_{S^1} \, f(\ps + \rM \omega)\, ds(\omega)\, d\rM \right)
    d\rT \, ds(\ps)  \\
    & =
    \frac{1}{2 \pi} \int_S \int_0^\infty
    \left( \int_0^\rT
    \int_{S^1}
    \frac{f(\ps + \rM \omega)}{\sqrt{\rT^2 - \rM^2}} \,  \rM \, d\rM \,
    dS(\omega) 
    \right) \, G(\ps,\rT) \, d\rT \, ds(\ps)
    \\
    & =
    \frac{1}{2 \pi} \int_S \int_0^\infty
    \left( \int_{\R^2}
    \frac{f(y)}{\sqrt{\rT^2 - \abs{y-p}^2}}\chi(\{\abs{y-p}<t\}) \, dy
    \right) \, G(\ps,\rT) \, d\rT \, ds(\ps)
    \\
    & =
    \frac{1}{2 \pi} \int_{\R^2} f(y)
    \left(
        \int_S \int_{\abs{y-\ps}}^\infty
            \frac{G(\ps,\rT)}{\sqrt{\rT^2 - \abs{y-\ps}^2}}
            \,d\rT \, ds(\ps)
    \right)
    \, f(y) \, dy
    \\ & =
    \langle f, \, \Po^\ast G \rangle,
\end{align*}
where
\beq{Gstar}
    \bigl(\Po^\ast G \bigr)(y)
    :=
    \frac{1}{2 \pi} \int_S \int_{\abs{y-\ps}}^\infty
    \frac{G(\ps,\rT)}{\sqrt{\rT^2 - |y-\ps|^2}} \, d\rT\, ds(\ps) \,.
\eeq
The integral in \req{Gstar} will be absolutely convergent for continuous $G$ 
provided that $G$ has a small amount of decay as $\rT \to \infty$, for
example if 
$G(\ps, \rT) = \mathcal{O}(1/\rT^{\alpha})$, as $\rT \to \infty$,  for some 
$\alpha >0$.

Next, we note a differentiation formula for the fractional integral
appearing in (\ref{ndsol}).

\begin{pr} Let $h$ be differentiable on $[0,\infty)$. 
Then, for $\rT>0$,
\beq{diffAbel}
    \di_\rT 
    	\int_0^{\rT} \frac{\rM \, h(\rM)}{\sqrt{\rT^2 - \rM^2}}  \,d\rM
    =
    \frac{1}{\rT} \int_0^{\rT}
    \frac{ \rM \,  (\di_\rM \ro h)(\rM)}{\sqrt{\rT^2-\rM^2}}  \,d\rM.
\eeq
\end{pr}

\begin{proof} 
Making the change of variable $\rM = \rT \xi$ in the integral 
on the left we have to evaluate
\[
\di_{\rT}     
	\int_0^1 \frac{\xi}{\sqrt{1-\xi^2}}
     	\, \rT h(\rT \xi)  \, d\xi.
\]
 Here differentiation under the integral yields
$
\int_0^1 \frac{\xi}{\sqrt{1-\xi^2}}
 \left( \rT \xi h'(\rT \xi) +h(\rT \xi) \right) 
\,d\xi,$ 
which is equal to the expression on the right side after changing back
to integration with respect to $\rM = \rT\xi$.
\end{proof}

\begin{proof5}
We compute $(\Po^\ast \tio \di_{\rT}^2  \Po f)(x)$ for smooth $f$ supported in
$\overline B$ and $x\in B$.  
The function $\tio \di_{\rT}^2 \Po f$ has decay of order
$1/ \rT^{2}$ as $\rT \to \infty$, and so lies in the domain of $\Po^\ast$.  
Using the definitions of $\Po$ and $\Po^\ast$, and relation
\req{diffAbel},
\begin{align*}
    (\Po^\ast \rT &\di_{\rT}^2 \Po f)(x)
    \\ 
    &=
    \frac{1}{2 \pi}
    \int_S \int_{\abs{x-\ps}}^\infty
    \di_{\rT}^2 \,
    \left( 
    	\int_0^{\rT} 
	\frac{\rM \, (\M f)(\ps,\rM)}{\sqrt{\rT^2-\rM^2}}
        \, d\rM
      \right)
     \frac{\rT \, d\rT \,ds(\ps) }{\sqrt{\rT^2 - \abs{x-\ps}^2}}
    \\ 
    &=
    \frac{1}{2 \pi}
    \int_S \int_{\abs{x-\ps}}^\infty
    \di_{\rT}
    \left( 
    	\frac{1}{\rT} 
    	\int_0^\rT 
	\frac{\rM \, (\di_{\rM} \ro \M f)(\ps,\rM)}{\sqrt{\rT^2-\rM^2}} 
     d\rM
     \right) \,
     \frac{\rT\, d\rT \,ds(\ps) }{\sqrt{\rT^2 - \abs{x-\ps}^2}}.
\end{align*}
Carrying out the differentiation in $\rT$ using the chain rule, using  
again \req{diffAbel}, and combining terms, the last integral can be
rewritten as 
\[
    \frac{1}{2\pi}
    \int_S 
    	\int_{\abs{x-\ps}}^\infty
    	\left( 
	\int_0^\rM
    		\frac{ \rM \, (\di_{\rM} \ro \di_{\rM} \ro \M f)(\ps,\rM)  
		-  \rM \, (\di_{\rM} \ro \M f)(\ps,\rM) }
		{\rT \sqrt{\rT^2 - |x-\ps|^2}   \sqrt{\rT^2 - \rM^2}}
		\, d\rM
    \right)
     d\rT \,ds(\ps) .
\]
Using the identity
\[
\di_{\rM} \ro \di_{\rM} \ro h
-
\di_{\rM} \ro h
= 
\di_{\rM} \ro ( \di_{\rM}  \ro h - h)  
= 
\di_{\rM} \ro   \ro \di_{\rM} h 
= 
\di_{\rM} \ro^2  \di_{\rM} h 
\]
and applying Fubini's theorem $ (\Po^\ast \rT \di_{\rT}^2 \Po f)(x)$ is
in turn is equal to
\[
    \frac{1}{2 \pi} 
    \int_S
    \int_0^\infty 
    \rM \,  ( \di_{\rM} \ro^2  \di_{\rM} \M  f)(\ps, \rM) 
    \left(
    \int_{\max(|x-\ps|,\rM)}^\infty
    \frac{d\rT}{\rT \, \sqrt{\rT^2 - |x-\ps|^2} \, \sqrt{\rT^2 - \rM^2}}
    \right)
    d\rM \,ds(\ps).
\]
The inner integral  evaluates to
\[
    \frac{1}{2 \rM \abs{ x-\ps}}
    \log\frac{\rM + \abs{x-\ps}}{ \wabs{\rM - \abs{x-\ps} }}
\]
giving
\begin{equation}\label{inner}
    (\Po^\ast \tio \di_{\rT}^2 \Po f)(x)
    =
    \frac{1}{4 \pi} \int_S
    \left(
    \int_0^\infty
    ( \di_{\rM} \ro^2  \di_{\rM} \M  f)(\ps, \rM) 
    \log \frac{\rM + |x-\ps|}{ \wabs{\rM - |x-\ps|} } d\rM
    \right) \frac{ds(\ps)}{\abs{x-\ps}} .
\end{equation}
Treating the inner integral in principal value sense, and integrating by
parts, it is equal to the limit as $\eps \to 0$ of boundary terms
\[
    \left[  
    (\ro^2  \di_{\rM} \M  f)(\ps, \rM) 
    \log \frac{\rM + \abs{x-\ps}}{ \abs{x-\ps} -\rM}
    \right]_0^{\abs{x-\ps} -\eps} 
    +
    \left[
    (  \ro^2  \di_{\rM} \M  f)(\ps, \rM) 
    \log
    \frac{\rM +\abs{x-\ps} }{\rM-\abs{x-\ps} }
    \right]_{\abs{x-y} +\eps}^{\infty}
\]
plus the term
\[
I_\epsilon := 
-\int_{\R^+ \setminus [\abs{x-\ps} -\eps, \abs{x-\ps} +\eps]}
	(\ro  \di_{\rM} \M  f)(\ps, \rM) \,
        \rM \di_{\rM} \log \frac{\rM +\abs{x-y}}{\wabs{\rM -\abs{x-\ps} }} 
        \,d\rM.
\]
Using that $\M f$ is smooth, flat at $\rM=0$,  and of bounded support in
$(0,\infty)$, the limit of the boundary terms is zero. Using the identity
\[
    \rM \di_{\rM} \log \frac{\rM+\vert
    x-\ps\vert}{\vert \rM -\abs{x-\ps} \vert} = -\abs{x-\ps}
    \di_{\rM}\log \vert \rM^2 -\abs{x-\ps}^2\vert
\]
followed by another integration by parts, yields the sum of another pair
of boundary terms and
    \[
        I_\epsilon  =
        -\abs{x-\ps}
        \int_{R^+\setminus [\abs{x-p} -\eps,\vert x-p \vert +\eps]}
        \left( 
        \di_{\rM} \ro \di_{\rM} \M f
        \right)(\ps,\rM)
        \log \wabs{ \rM^2 -\abs{x-\ps}^2}
    \,d\rM.
    \]
    The boundary terms again evaluate to zero as $\eps \to 0$ while
    the integral $I_\epsilon$ converges to
    \[
    -\abs{x-\ps} \int_0^{\infty}
    \left( 
        \di_{\rM} \ro \di_{\rM} \M f
        \right)(\ps,\rM)
 \log \wabs{\rM^2-\abs{x-\ps}^2 } \,d\rM.
    \]
    Inserting this into (\ref{inner}) and taking into account the support of
    $\M f$, gives
\[
	(\Po^\ast \rT \di_{\rT}^2 \Po f )(x)
   	= -\frac{1}{4\pi}
    	\int_S \int_0^{2\Rball}
      	\left( 
        		\di_{\rM} \ro \di_{\rM} \M f
        \right)(\ps,\rM)
        \log \wabs{ \rM^2-\abs{x-\ps}^2 } \, d\rM\,ds(\ps).
\]
In view of \req{m-inv} of Theorem \ref{thm:circinv}, 
\req{invP} in Theorem \ref{thm:invPW} is proved, for $n=2.$ \eproof
\end{proof5}

\begin{proof6} 
	Formula (\ref{asymm}), for $n=2$, 
	an easy corollary of the result just established. Indeed, for $f,g$
	smooth 
	with compact support in the closed disk of radius $\Rball$, then
	\begin{equation}\label{trace}
		\langle f, g\rangle
		=
		- \frac{2}{\Rball}
		\left\langle
			\Po^\ast \rT \di_{\rT}^2 \Po f, \, g
		\right\rangle
		=
		- \frac{2}{\Rball}
		\left\langle
			\rT \di_{\rT}^2 \Po f, \, \Po g
		\right\rangle \,,
	\end{equation}
	which is (\ref{asymm}) for $n=2$, due to the definition of the operator
	$\Po$.

	In (\ref{asymm}), the left hand side is
	symmetric in $f$ and $g$, while the right side is not. Thus there is a
	companion identity, reversing the roles of $u$ and $v$ on the
	right. Taking the difference gives the equation
	\[ 
	0=
	\int_S \int_0^{\infty} 
	\rT (u_{\rT\rT} v - u v_{\rT\rT}) \, d\rT \, ds(\ps).
	\] 
	Integrating by parts (the boundary terms vanish) yields
	\[
	0 =
	\int_S \int_0^{\infty}
	\left( u_{\rT} v - u v_{\rT} \right)
	\, d\rT \, ds(\ps),
	\]
	and another integration by parts proves
	\[
	0   =
	\int_S \int_0^{\infty}
	u_{\rT} v
	\,d\rT\,ds(\ps)
	= \int_S\int_0^{\infty}
	u v_\rT \,d\rT\,ds(\ps).
	\]
	Using this and one integration by parts in (\ref{asymm}) establishes
	(\ref{symm}), which completes the proof of Theorem \ref{thm:e2m} for
	$n=2$. The extension to higher (even) dimensions follows almost word
	for word the proof from \cite[Section 4.2]{FinPatRak04}, where the
	trace identities in odd dimensions greater than three were proved
	from the three dimensional case.  
 \eproof
\end{proof6}

\begin{proof7}
Reversing the chain of reasoning in (\ref{trace}) proves \req{invP} in $L^2$
sense from
(\ref{asymm}). Similarly, \req{invW} follows from (\ref{symm}). 
However, as both sides are continuous functions when $f$ is smooth, 
the formulas hold pointwise as well. 
\eproof
\end{proof7}

\section{Numerical results}

In the previous sections, we have established several
exact inversion formulas to recover a function $f$
supported in a closed disc $\overline D$ from either its spherical
means $\M f$ or the trace $\W f$ of the solution of the wave equation
with initial data $(f,0)$.  However, those formulas require continuous
data, whereas in practical applications only a discrete data set is
available.  For example, in thermoacoustic tomography (see Figure
\ref{fg:tact-line}) only a finite number of positions of the line
detectors and finite number of samples in time are feasible. In this
section we derive discrete \textit{filtered back-projection} (FBP)
algorithms with linear interpolation in dimension two and present some
numerical results.

The derived FBP algorithms are numerical implementations of
discretized versions of \req{m-lap}-\req{filbac} and
\req{wavefinite}-\req{invW} and the derivation of any of them follows
the same line.  We shall focus on the implementation of \req{m-inv},
assuming uniformly sampled discrete data
\beq{md}
     F^{k,m}  &:= (\M f ) ( \ps^k , \rM^m  ) \,, \quad &(k,m) \in   \set{ 0,
       \dots,  N_\ph } \times \set{ 0, \dots,  N_{\rM}}\,, 
\eeq
where $\ps^k :=  \Rball \left( \cos (k h_{\ph})  ,   \sin (k h_{\ph})
\right) $, 
$\rM^m := m h_{\rM}$, $ h_{\ph} := 2 \pi / (N_{\ph}+1) $  and 
$h_{\rM} := 2\Rball/N_{\rM}$.
In order to motivate the  derivation of a discrete FBP algorithm  based on
\req{m-inv}, 
we introduce the   differential operator $\Do := \di_{\rM} \ro  \di_{\rM}$
and the integral operator 
\beq{Log}
    \Io : C_0^\infty (  S \times [0, 2\Rball) ) \to C^\infty (  S \times [0,
	2\Rball) ) 
    \\
    (\Io G)(\ps, \rH)
    := \int_{0}^{2\Rball} G(\ps, \rM)  \log  \abs{ \rM^2 - \rH^2}  d\rM
\eeq
which both act in the  second component,  and the so called
\textit{back-projection operator}
\beq{bp}
    \B : C^\infty (  S \times [0, 2\Rball) ) &\to C^\infty( \overline D)
    \\
    (\B G)(x)
    &:= \frac{1}{2\pi\Rball} \int_{S} G(\ps, \abs{x-\ps}) ds(\ps)
    \\
    &= \frac{1}{2\pi} \int_{0}^{2\pi} G( \ps(\ph), \abs{x - \ps(\ph)}) d\ph
    \;,
\eeq
where $\ps(\ph) := \Rball (\cos\ph, \sin\ph)$.  Therefore, we can rewrite
\req{m-inv} as 
\beq{m-inv1}
    f = ( \B  \Io  \Do) (\M f)  \,.
\eeq
In the numerical implementation the operators $\B$, $\Io$,  and $\Do$ in
\req{m-inv1} 
are replaced with  finite dimensional approximation $\Bd$,  $\Id$ and $\Dd$
(as described below) and 
\req{m-inv1} is approximated by
\beq{m-inv-d}
    f(x^i) \approx f^i := ( \Bd  \Id  \Dd  \f F )^i \,, \quad i \in \set{0,
      \dots, N}^2 \,. 
\eeq
Here $\f F := (F^{k,m})_{k,m}$ with $F^{k,m}$ defined by \req{md}, $x^{i} :=
- (\Rball,\Rball) + i h_x$ with $i =(i_1, i_2) \in \set{0, N}^2$ and $h_x :=
2\Rball/N$.  In the following $\Sd_{\ph,r}$ and $\Sd_{x}$ denote the
sampling operators that map $G \in C^{\infty}(S\times [0, 2\Rball])$ and $f
\in C^{\infty}(\overline D)$ onto its samples, $\Sd_{\ph,r} G : = (G(\ps^k,
\rM^{m}))_{k,m}$ and $\Sd_{x} f := \f f := (f(x^i))_{i}$, where we set
$f(x^i) := 0$ if $x^i \not \in \overline D$.  Moreover,
$\abs{\cdot}_{\infty}$ denotes the maximum norm on either $\R^{(N_\ph+1)
\times (N_{\rM}+1)}$ or $\R^{(N+1) \times (N+1)}$.

\begin{enumerate}
\item
The operator  $\Do$ can be written as  $ \di_{\rM} +  \ro \di_{\rM}^2$.
We approximate $\di_{\rM} G$ with symmetric  finite differences
$\bigl( G^{k,m+1} - G^{k,m-1} \bigr) / (2h_{\rM})$,
$\di_{\rM}^2 G $ by $\bigl( G^{k,m+1} + G^{k,m-1}  - 2 G^{k,m}\bigr) / h_\rM^2$
and  the multiplication operator $G \mapsto \rM  G $ by point-wise discrete 
multiplication $(G^{k,m})_{k,m} \mapsto (\rM^m G^{k,m})_{k,m}$.
This leads to the discrete approximation
\beq{filt1}
    &\Dd:
    \R^{(N_\ph+1) \times (N_{\rM}+1)}
    \to
    \R^{(N_\ph+1) \times (N_{\rM}+1)}\,, \\
    &\f F \mapsto
    \left(
        (\Dd \f G) ^{k,m}  :=
        \frac{1}{h_{\rM}}
        \left(
            \bigl( m +\frac{1}{2} \bigr) G^{k, m+1} +
            \bigl( m -\frac{1}{2} \bigr)  G^{k, m-1} -
            2m G^{k,m}
        \right)
    \right)_{k,m}
\eeq
where  we set $G^{k,-1} := G^{k,N_{\rM}+1} := 0$.
The approximation of $\di_{\rM}$ with symmetric  finite differences is  of
second 
order and therefore $\abs{ (\Sd_{\ph,r} \Do - \Dd \Sd_{\ph,r}) G}_{\infty}
\leq  C_1 h_\rM^2$ for some constant $C_1$  which does not depend on
$h_{\rM}$. 

\item
Next we  define a second order approximation to the integral
operator $\Io$.   This is done by replacing $G(\ps^k, \cdot)$  in \req{Log}
by the 
piecewise linear  spline $T^k[G]:  [0, 2\Rball] \to \R$ interpolating $G$ at
the  nodes 
$\rM^m$. More precisely,
\beqn
    \Id:
    \R^{(N_\ph+1) \times (N_{\rM}+1)}
    \to
    \R^{(N_\ph+1) \times (N_{\rM}+1)}:
    \quad
    \f G \mapsto
    \left(
        (\Id \f G) ^{k,m}
    \right)_{k,m}
\eeqn
is defined by
\beq{linSpline}
    T^k[\f G](\rM)
    :=
    G^{k,m} + \frac{\rM-\rM^m}{h_{\rM}} ( G^{k, m+1} - G^{k, m} ) \;,
    \quad
    \rM \in [\rM^{m}, \rM^{m+1}]
\eeq
and
\beq{LogD}
    (\Id \f F)^{k,m}
        & :=
        \int_{0}^{2 \Rball}
            T^k[\f G](\rM)
            \log \abs{ \rM^2  - (\rM^{m})^2}
        d\rM
        \\ & =
        \sum_{m'=0}^{N_{\rM}-1}
        G^{k,m'}
        \left(
            \int_{\rM^{m'}}^{\rM^{m'+1}}
                \log \abs{\rM^2 - (\rM^{m})^2 }
            d\rM
        \right)
        \\ &+
        \sum_{m'=0}^{N_{\rM}-1}
        \frac{G^{k, m'+1}- G^{k, m'}}{h_{\rM}}
        \left(
            \int_{\rM^{m'}}^{\rM^{m'+1}}
                (\rM - \rM^{m'})\log \abs{\rM^2 - (\rM^{m})^2}
            d\rM
        \right)
    \;.
\eeq
For an efficient  and  accurate numerical implementation it is crucial  that
the integrals in \req{LogD} are evaluated analytically.   In fact,
by straight forward  computation it can be verified that
\beq{LogD2}
	(\Id \f F)^{k,m}
	&
	=
	\sum_{m'=0}^{N_{\rM}-1}
		a^{m}_{m'} G^{k,m'}
		+
		\frac{1}{h_{\rM}}
		\sum_{m'=0}^{N_{\rM}-1}
			b^{m}_{m'}
			\left(
				G^{k, m'+1}- G^{k, m'}
			\right)
	\;,
	\\
	a^{m}_{m'}
	&:=
	\Bigl[
		( \rM - \rM^{m} )
		\log \abs{ \rM - \rM^{m} }
		+
		(\rM + \rM^{m})
		\log \abs{\rM + \rM^{m}}
		-
		2 \rM
	\Bigr]_{\rM=\rM^{m'}}^{\rM^{m'+1}}
	\,,
	\\
	b^{m}_{m'}
	&:=
	- \rM^{m'} a^{m}_{m'}
	+
	\frac{1}{2}
	\Bigl[
		(\rM^2 - (\rM^{m})^2)
		\log \abs {\rM^2 - (\rM^{m})^2}
		-
		\rM^2
	\Bigr]_{r=\rM^{m'}}^{\rM^{m'+1}}
	\,.
\eeq
Moreover, using the fact that piecewise linear interpolation is of second
order \cite{QSS00} and that $ \rM \mapsto \log  \abs{ \rM - (\rM^{m})^2}$ is
integrable, 
it can be readily verified that  the approximation error satisfies
$\abs{ (\Sd_{\ph,r} \Io  - \Id \Sd_{\ph,r} ) G}_{\infty} \leq  C_2 h_\rM^2$
with some constant $C_2$ independent of $h_{\rM}$.

\item
Finally, we define a second order approximation to the
back-projection \req{bp}. 
The discrete back-projection  operator 
$\Bd: \R^{(N_\ph+1) \times (N_{\rM}+1)}   \to \R^{(N+1) \times (N+1)}$
is  obtained by approximating \req{bp} with the trapezoidal rule
and piecewise linear interpolation \req{linSpline} in the second
variable,
\beq{bpD}
	(\Bd \f G)^i
	:=
	\frac{1}{N_\ph +1}
	\sum_{k=0}^{N_\ph}
		T^k[\f G] ( \abs{x^i - \ps^k} ) 
		\,,
		\quad x^i \in D,
\eeq
and setting $(\Bd \f G)^i := 0$ for $x^i \not \in D$.
It is well known \cite{QSS00} that both linear interpolation in $r$ and the
trapezoidal rule  
in $\ph$ are second order approximations and therefore 
$\abs{ (\Sd_{x} \B  - \Bd \Sd_{\ph,r}) G}_{\infty} \leq  C_3 \max \set{
  h_\rM^2,  h_\ph^2}$ 
for some constant $C_3$.
\end{enumerate}

\begin{algorithm}[H]
\caption{Discrete FBP algorithm with linear interpolation
for reconstruction  $\f f $ using  data $\f F$. }
\label{alg}

\begin{alginc}

\State
$ h_\ph \gets 2\pi/(N_\ph +1)$
\State
$ h_{\rM} \gets 2\Rball/N_{\rM}$  \Comment{initialization}

\For{$ m , m' = 0, \dots, N_{\rM}$} \Comment{Pre-compute kernel}
	\State
	{\tt Calculate $a_{m'}^m$, $b_{m'}^m$ according to \req{LogD}}
\EndFor

\State
\For{ $k =0,\dots,N_\ph$}  \Comment{Filtering}
	\For{ $m = 0,\dots, N_{\rM}$}
		\State $ F^{k,m} \gets
		\bigl( m + 1/2 \bigr) F^{k, m+1} +
		\bigl( m - 1/2 \bigr)  F^{k, m-1} -
		2m F^{k,m}  $ \Comment{Equation \req{filt1}}
	\EndFor
	
	\For{ $m = 0,\dots, N_{\rM}$}
		\State $ F^{k,m}
		\gets
		\sum_{m'=0}^{N_{\rM}-1}
		a^{m}_{m'} F^{k,m'}
		+
		\sum_{m'=0}^{N_{\rM}-1}
		b^{m}_{m'}
		\left(
			F^{k, m'+1}- F^{k, m'}
		\right)/h_{\rM}$
		\Comment{Equation \req{LogD}}
	\EndFor
\EndFor

\State

\For{ $i_1, i_2 = 0,\dots, N$}  \Comment{FBP with linear interpolation}
	\State 
        $i \gets (i_1,i_2)$
        \State
	$f^i  \gets 0$
	\For{ $k = 0, \dots, N_\ph$}
		\State
		{\tt Find $m  \in \set{0,\dots,N_{\rM}-1}$ with}
		$\rM^{m} \leq \abs{\ps^{k} - x^i} < \rM^{m+1} $
		
		\State
		$T \gets F^{k,m} + (\rM-\rM^{m})( F^{k, m+1} - F^{k, m}
		)/{h_{\rM}} $ 
		\Comment{interpolation \req{linSpline}}

		\State
		$f^i  \gets  f^i +  T/(N_\ph+1) $ \Comment{discrete
		back-projection \req{bpD}} 
	\EndFor 
\EndFor
\end{alginc}
\end{algorithm}

The discrete FBP algorithm is given by \req{m-inv-d} with $\Dd$, $\Id$,
$\Bd$ defined 
in \req{filt1}, \req{LogD2}, \req{bpD} and is summarized in Algorithm
\ref{alg}. 
Using 
	$f(x^i) = (\Sd_{x} \B  \Io  \Do F)^i =(\Sd_{x} f)^i$ 
and 
	$f^i = (\Bd  \Id  \Dd \Sd_{\ph,r} F)^i$,
the discretization error 
	$\abs{f(x^i)-f^i}$ 
can be estimated as
\beq{errDec}
	\abs{ (\Sd_{x} \B  \Io  \Do -   \Bd \Id  \Dd \Sd_{\ph,r}) F}_{\infty}
	& \leq
	\abs{ (\Sd_{x} \B - \Bd \Sd_{\ph,r}) ( \Io \Do F ) }_{\infty}
	\\
	&
	+
	\abs{ \Bd (\Sd_{\ph,r} \Io - \Id \Sd_{\ph,r})(\Do F) }_{\infty}
	\\
	&
	+
	\abs{ \Bd \Id (\Sd_{\ph,r} \Do - \Dd \Sd_{\ph,r})(F)}_{\infty}.
\eeq
Using the facts that $\Bd$ and $\Id$ are bounded by some  constant
independent of $h_{\rM}$ and 
that the approximation of  $\Do$, $\Io$, $\B$ with $\Dd$, $\Id$, $\Bd$ are
of second order, 
implies that
\beq{est}
	\abs{ \Sd_{x} f - \Bd \Id \Dd \f F }_{\infty}
	\leq
	C \max \set{ h_\rM^2, h_\ph^2 }
	\,,
\eeq
for some constant $C$ independent of $h_{\rM}$, $h_\ph$.
This shows that the derived FBP algorithm has  second order accuracy (for
exact data). 

\begin{figure}[htb!]
\begin{center}
\includegraphics[width = 0.49\textwidth, height =
  0.35\textwidth]{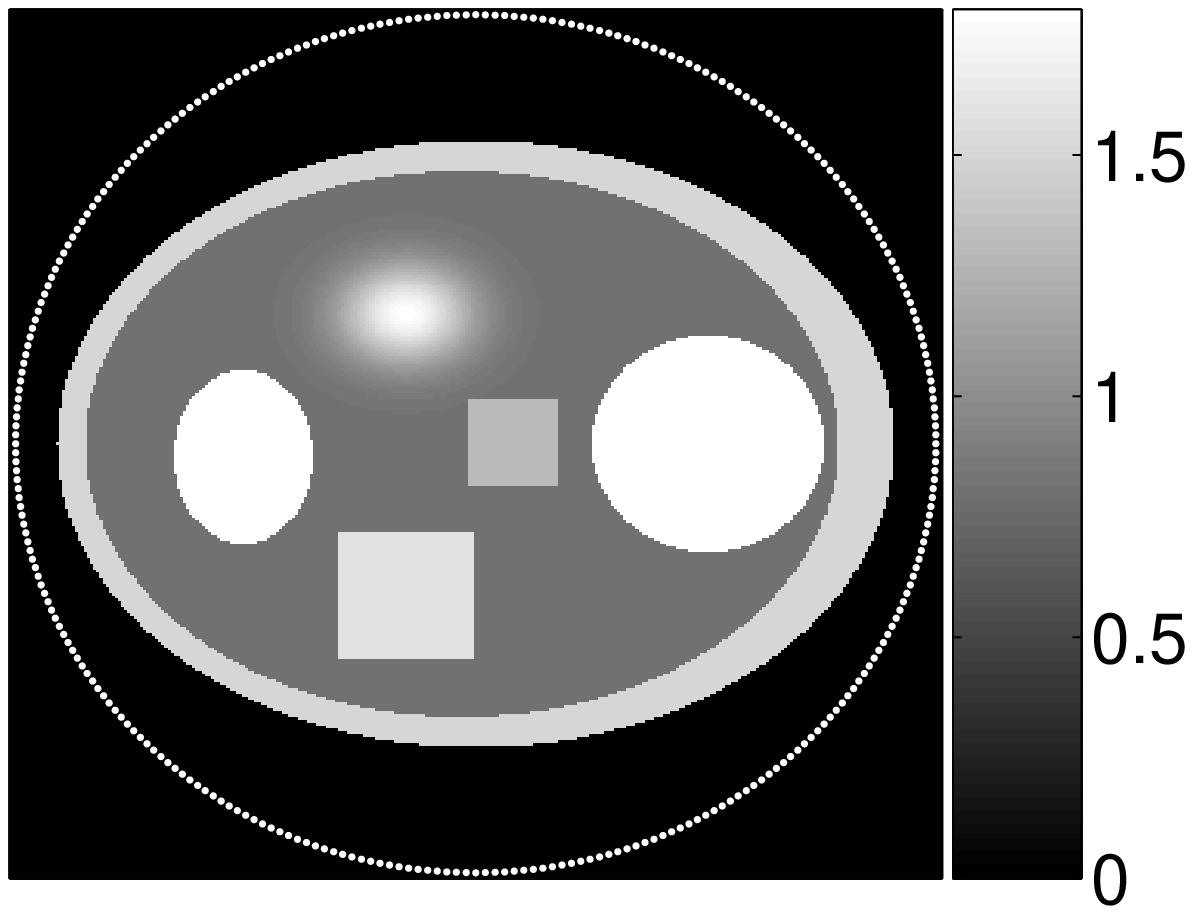} 
\includegraphics[width = 0.49\textwidth, height =
  0.35\textwidth]{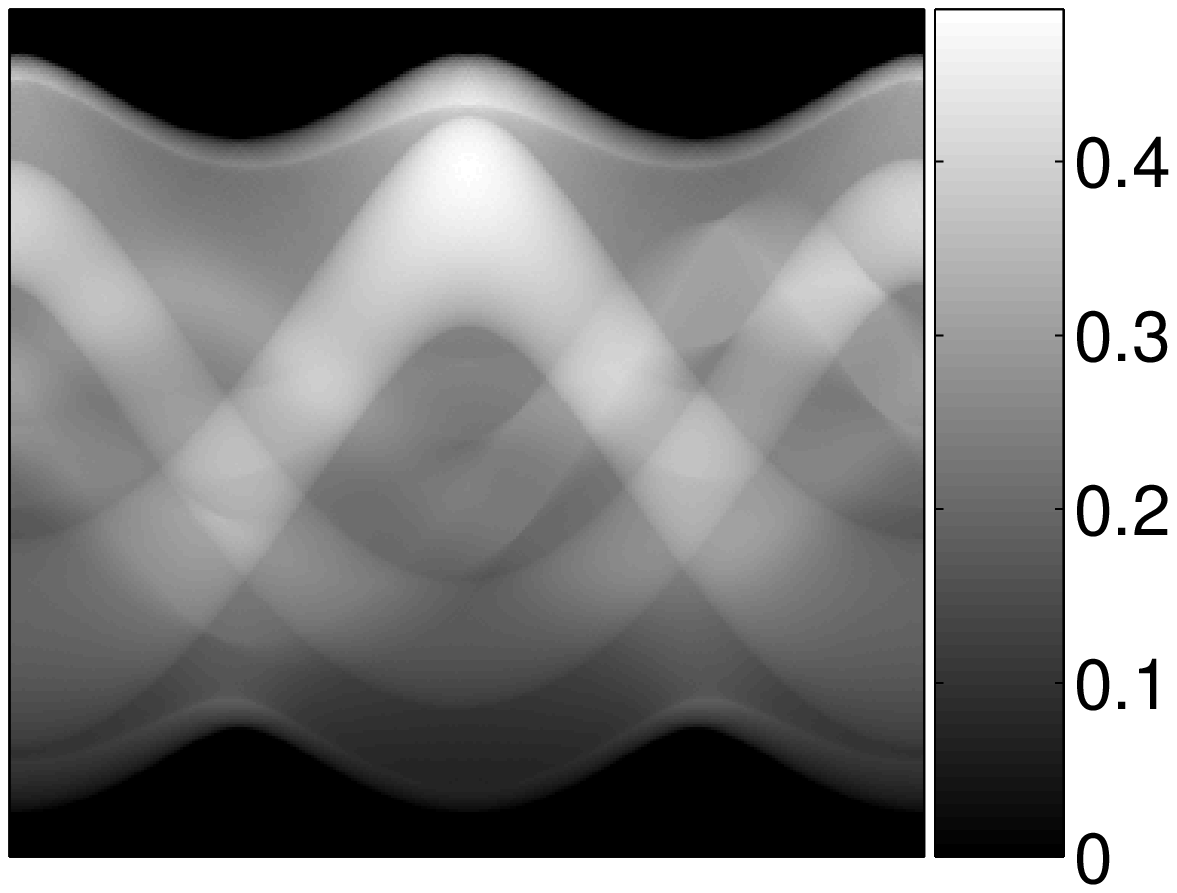} 
\end{center}
\caption{\textbf{Imaging phantom and data}.
Left:  Imaging phantom $f$ consisting of several characteristic functions
and one 
Gaussian kernel. Right: Simulated data $F = \M f$.}
\label{fig:phantom}
\end{figure}

In the numerical implementation, the coefficients in \req{filt1}, are
pre-computed and stored. Therefore the numerical effort of the evaluating
\req{filt1} is $\mathcal O(N_\rM^2 N_\ph)$.  Moreover, \req{filt1} requires
$\mathcal O(N_{\rM} N_\ph)$ operations and the discrete FBP $\mathcal O (N^2
N_\ph)$, since for all $(N+1)^2$ reconstruction points $x^i$ we have to sum
over $N_\ph+1$ center locations on $S$.  Hence, assuming $N \sim N_{\rM}$
and $N \sim N_\ph$, Algorithm 1 requires $\mathcal O (N^3)$ operations and
therefore has the same numerical effort as the classical FBP algorithm used
in x-ray CT \cite{Nat86}.  Analogous to the procedure described above,
discrete FBP algorithms were derived using equation \req{m-lap},
\req{m-inv-hilb} for inverting $\M$ and \req{wavefinite} for inverting $\W$.

\begin{figure}[htb!]
\begin{center}
\includegraphics[width = 0.5\textwidth, height =
  0.35\textwidth]{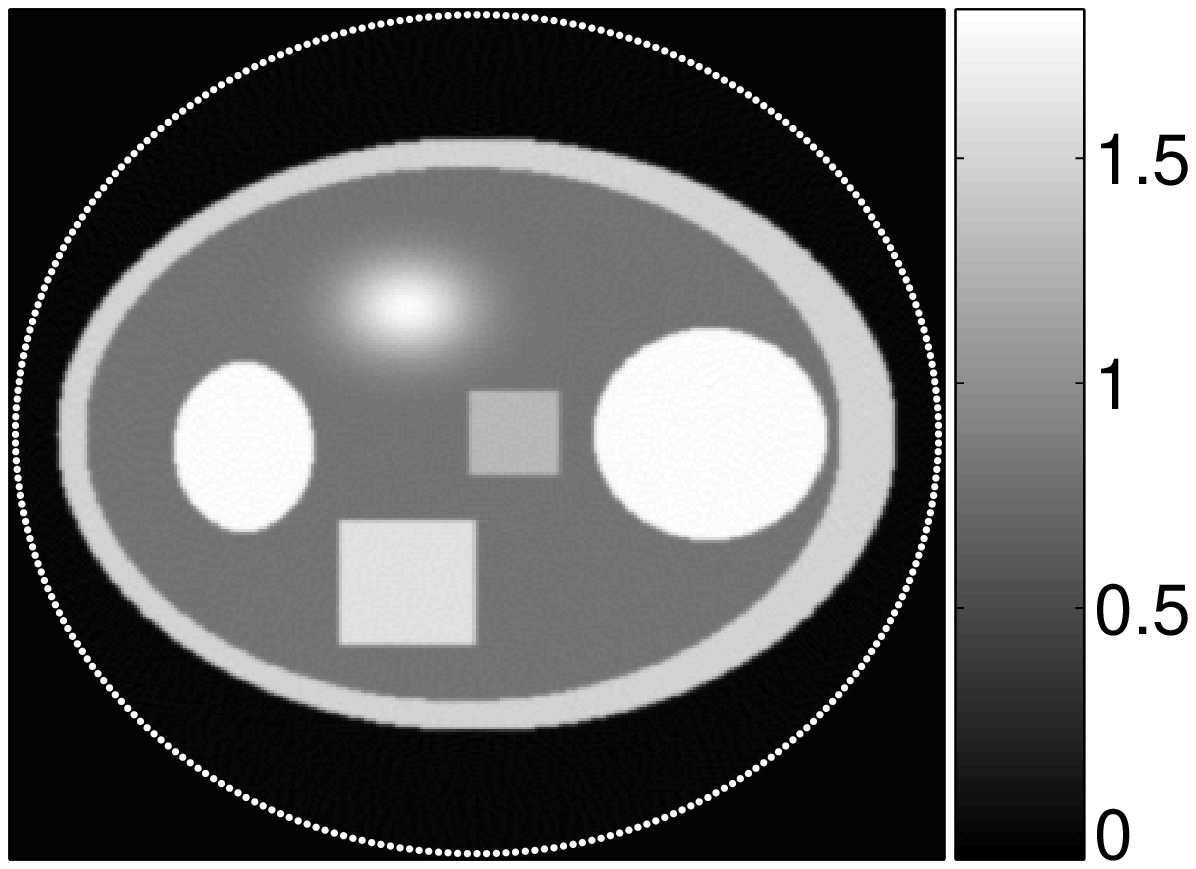} 
\includegraphics[width = 0.45\textwidth, height =
  0.35\textwidth]{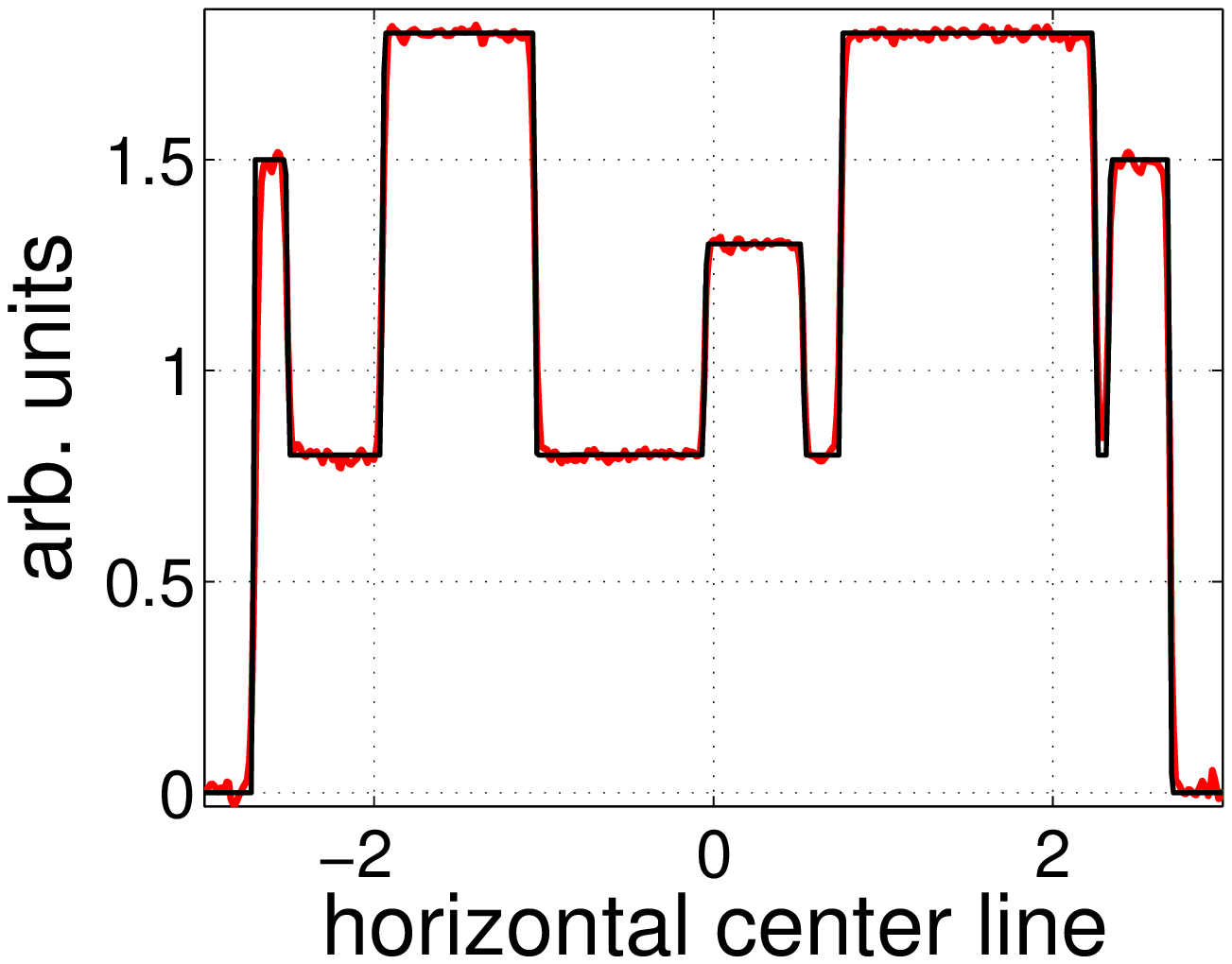}\\ 
\includegraphics[width = 0.5\textwidth, height =
  0.35\textwidth]{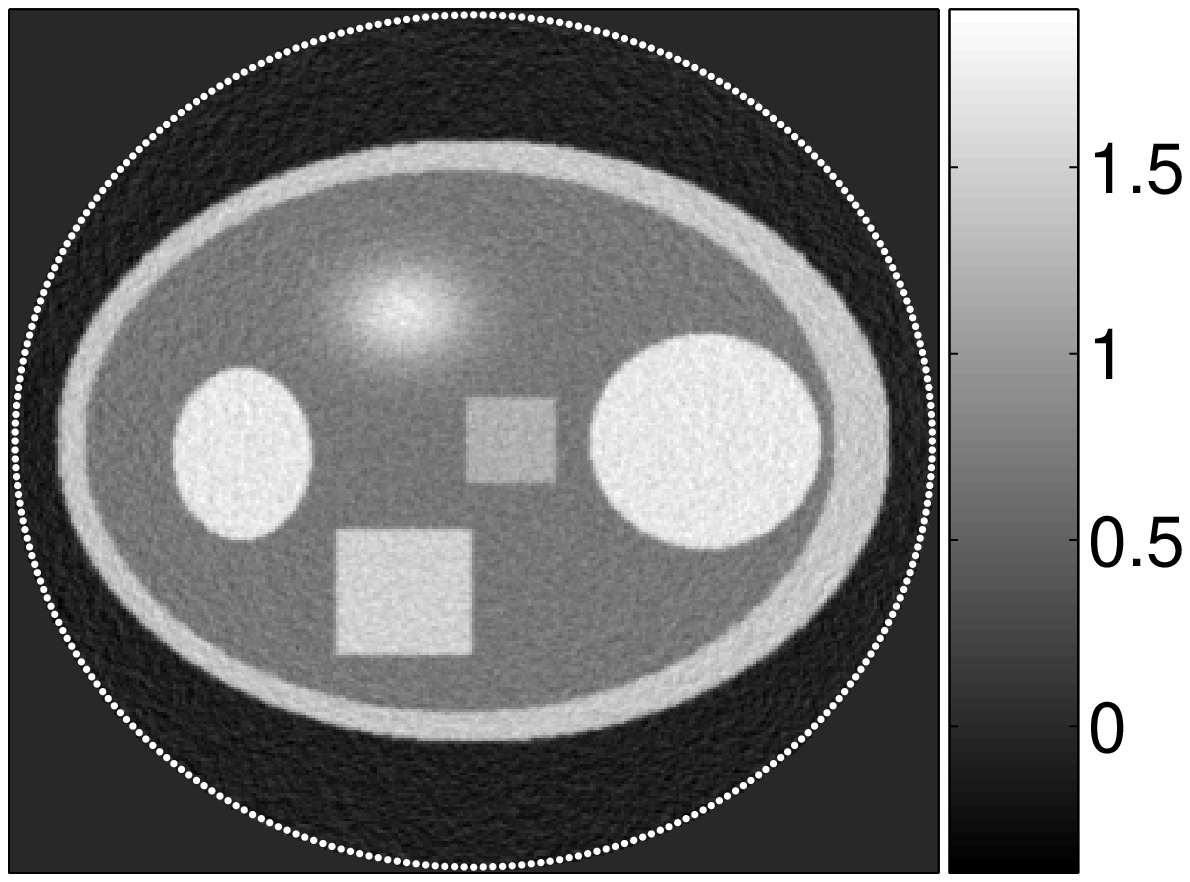} 
\includegraphics[width = 0.45\textwidth, height =
  0.35\textwidth]{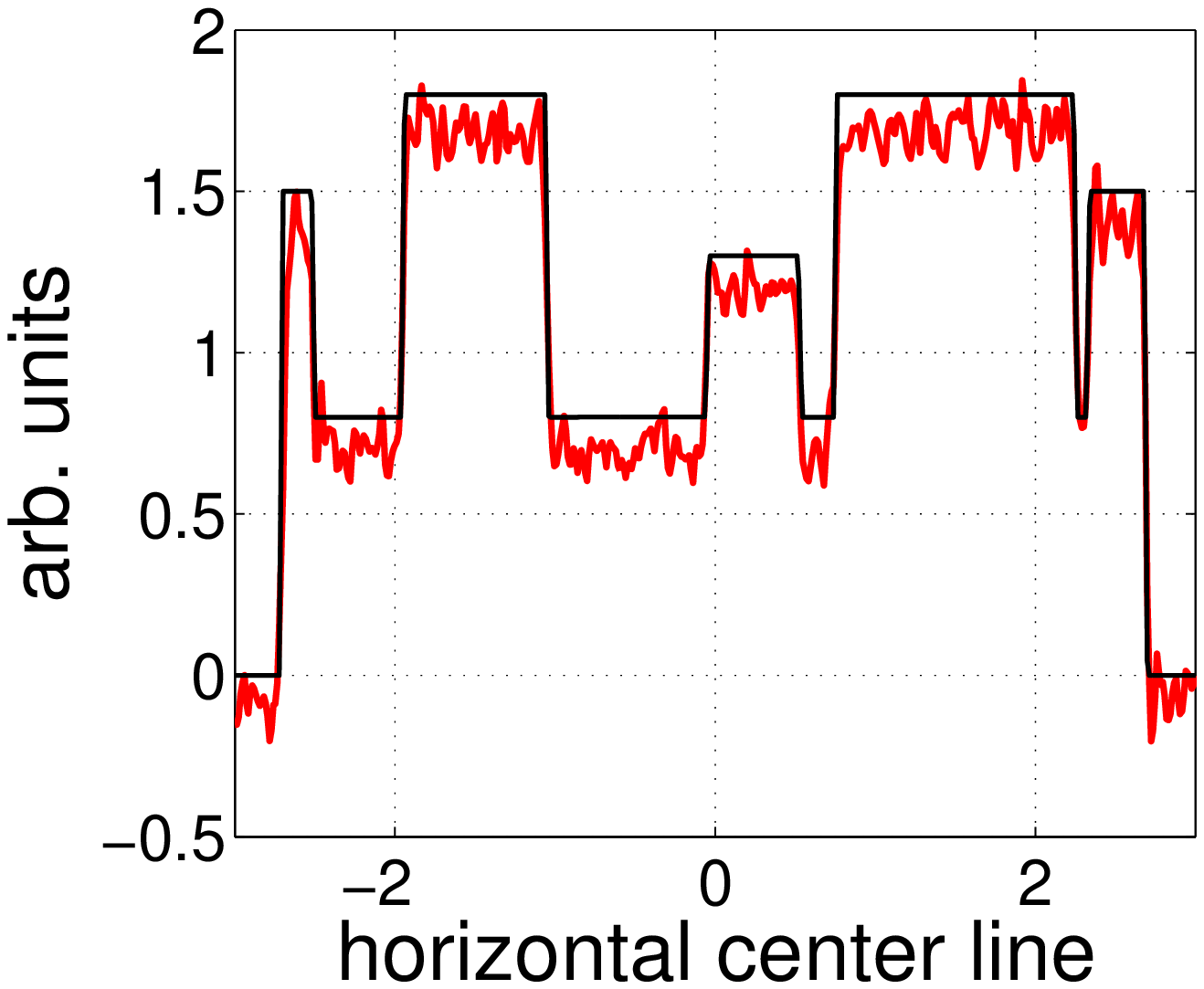}\\ 
\end{center}
\caption{\textbf{Numerical Reconstruction with Algorithm 1.}
Top: Reconstructions from simulated data. Bottom:  Reconstructions from
simulated data after adding $5 \%$ uniformly distributed noise.} 
\label{fig:log} 
\end{figure}

\begin{figure}[htb!]
\begin{center}
\includegraphics[width = 0.5\textwidth, height =
  0.35\textwidth]{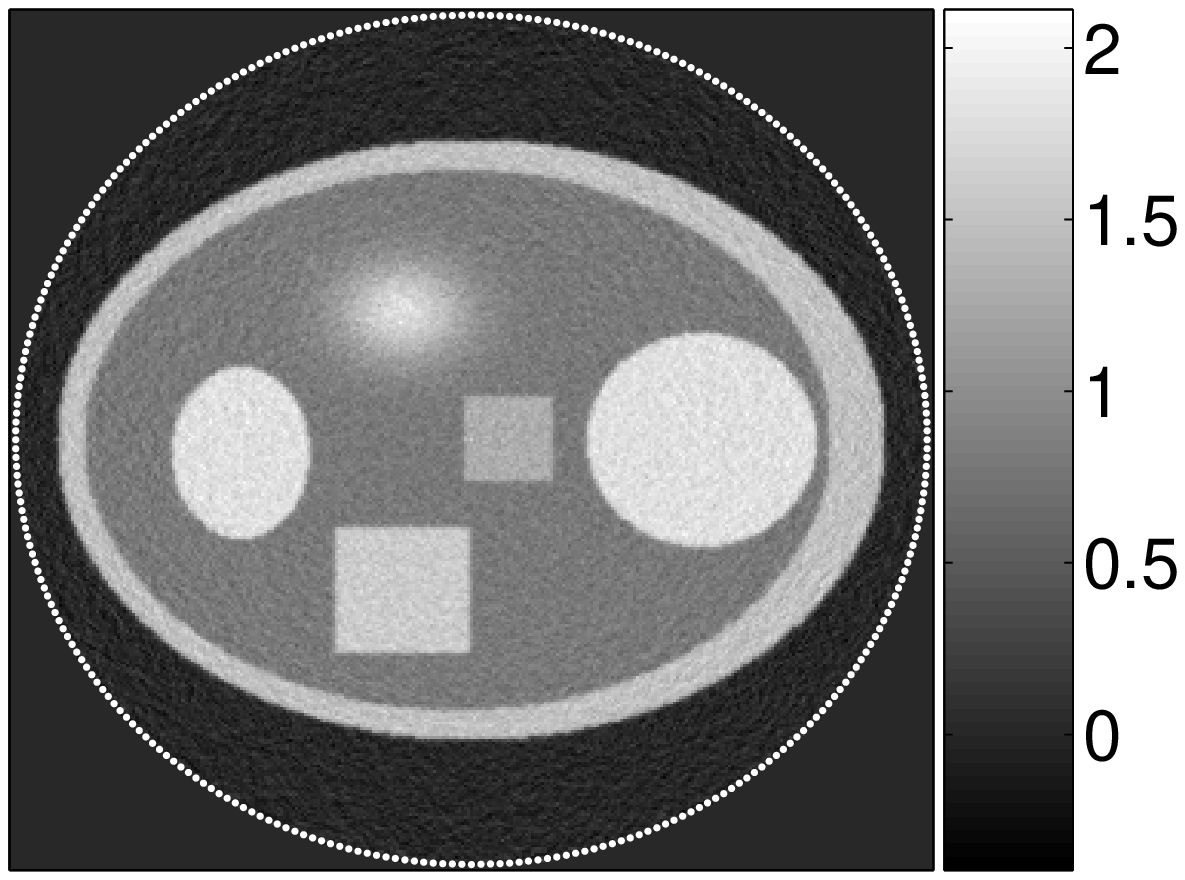} 
\includegraphics[width = 0.45\textwidth, height =
  0.35\textwidth]{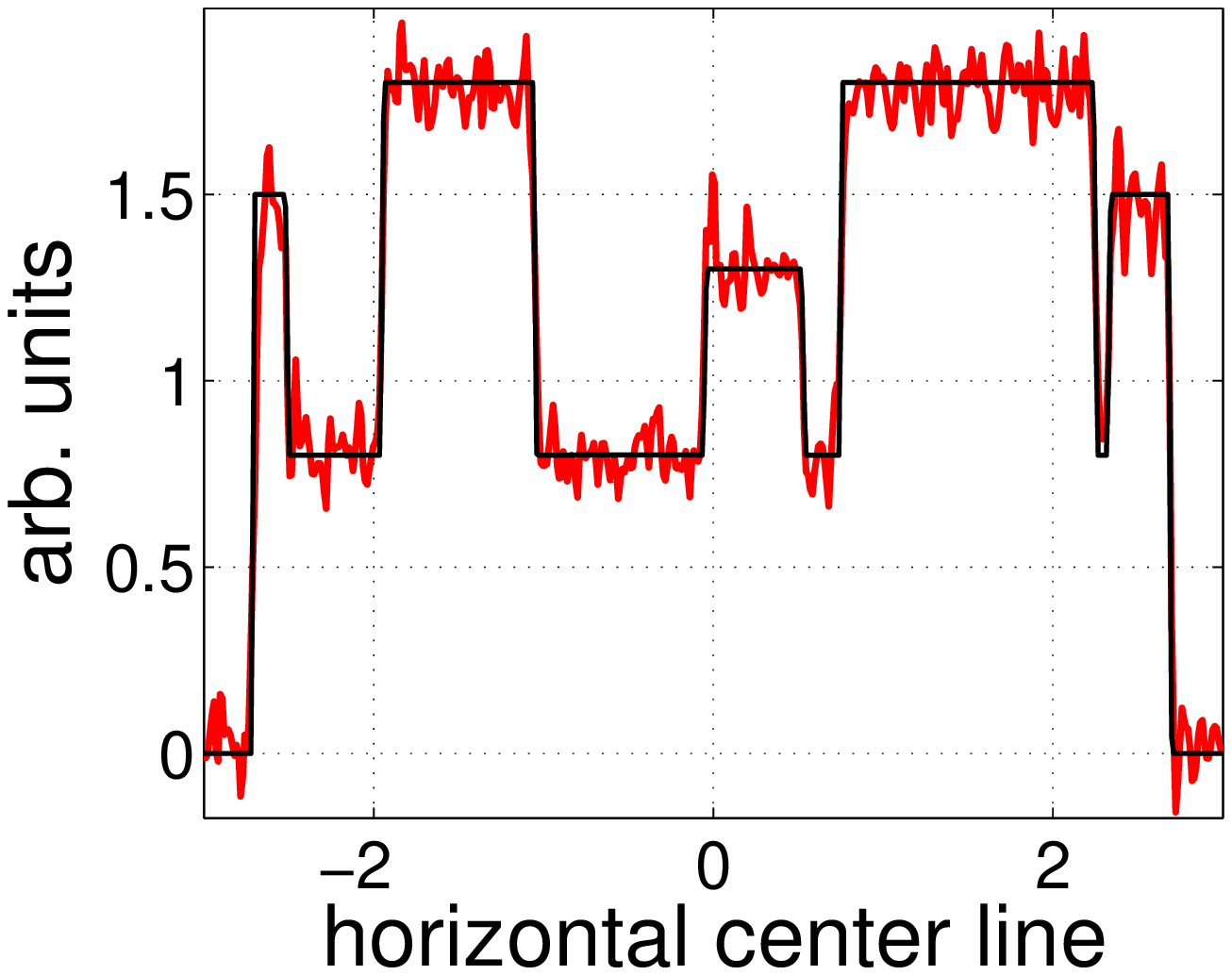}\\ 
\includegraphics[width = 0.5\textwidth, height =
  0.35\textwidth]{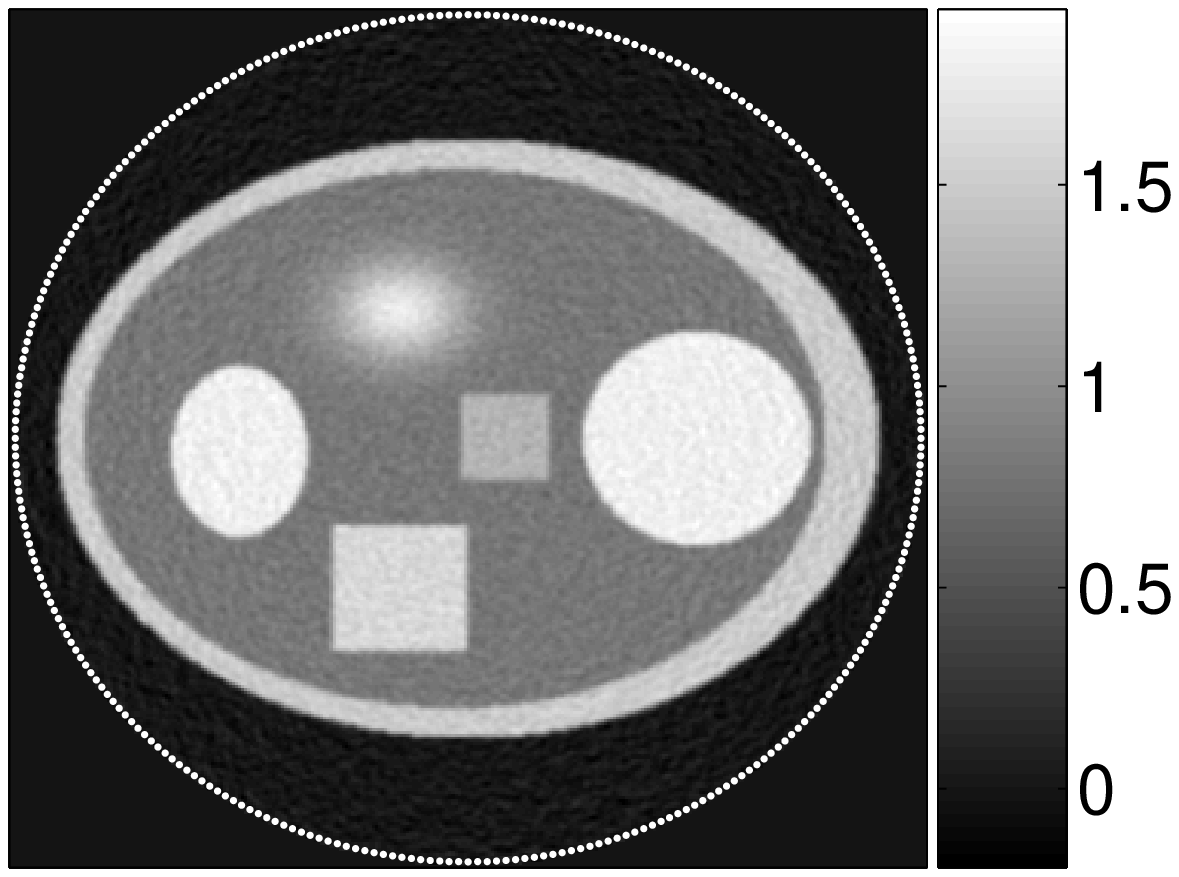} 
\includegraphics[width = 0.45\textwidth, height =
  0.35\textwidth]{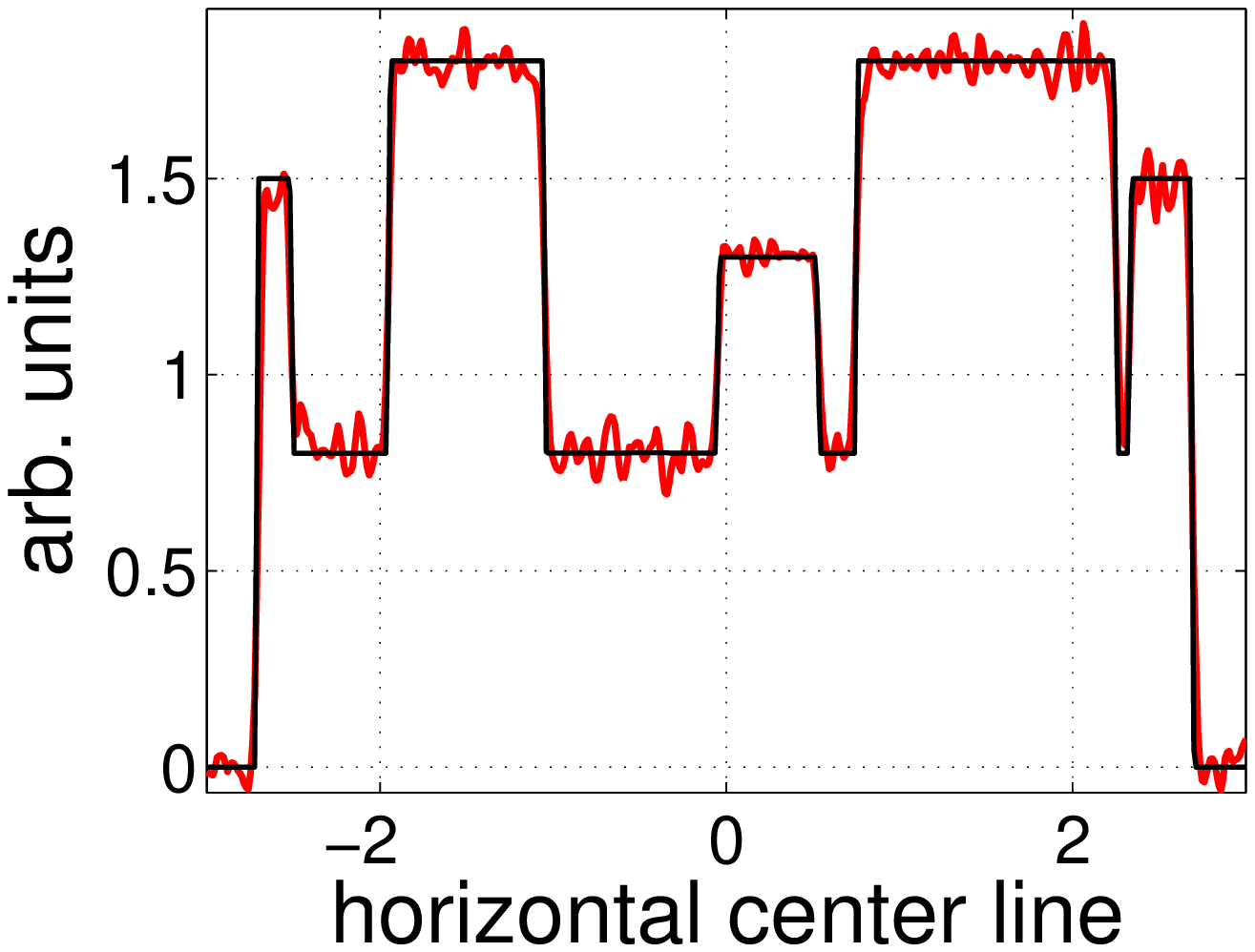}\\ 
\end{center}
\caption{\textbf{Numerical reconstruction  from spherical means with $5\%$
    noise added.} 
Top:  Reconstruction using  \req{m-lap}.  Bottom: Reconstruction using
    \req{m-inv-hilb}. } 
\label{fig:laplog}
\end{figure}

\begin{figure}[htb!]
\begin{center}
\includegraphics[width = 0.5\textwidth, height =
  0.35\textwidth]{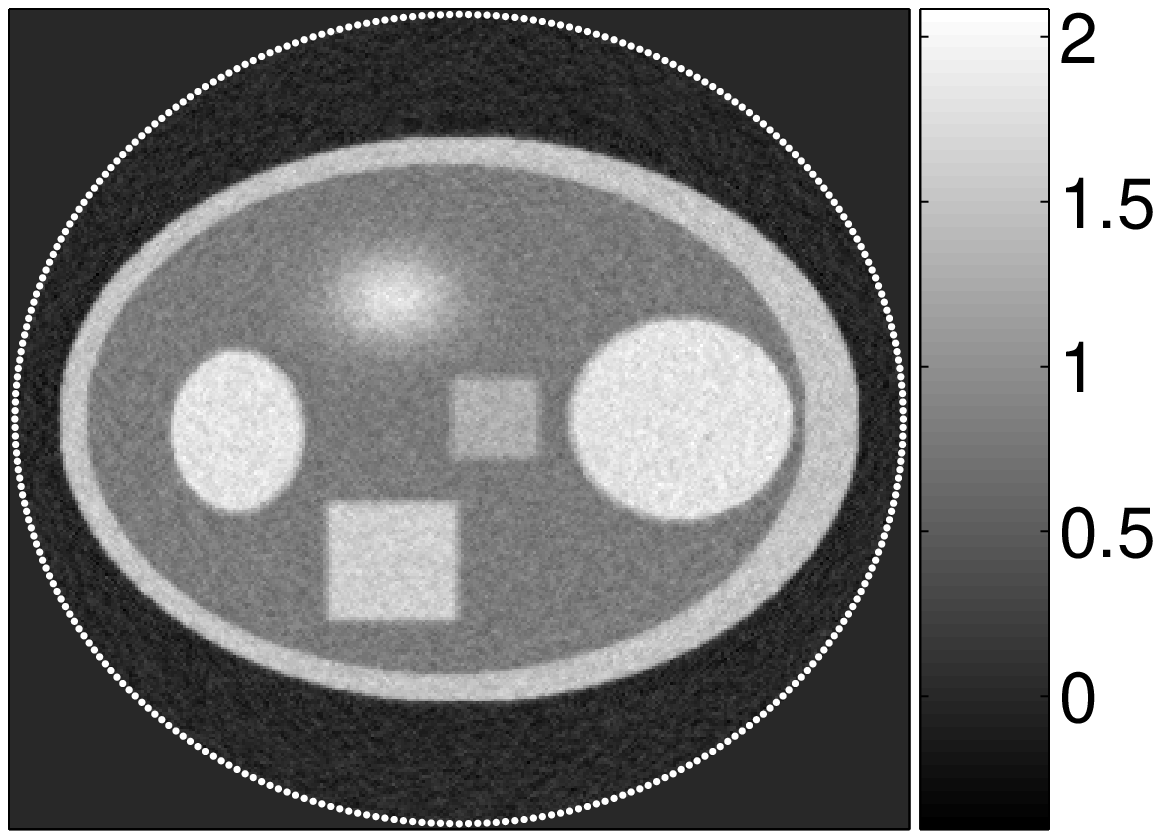} 
\includegraphics[width = 0.45\textwidth, height = 0.35\textwidth]
{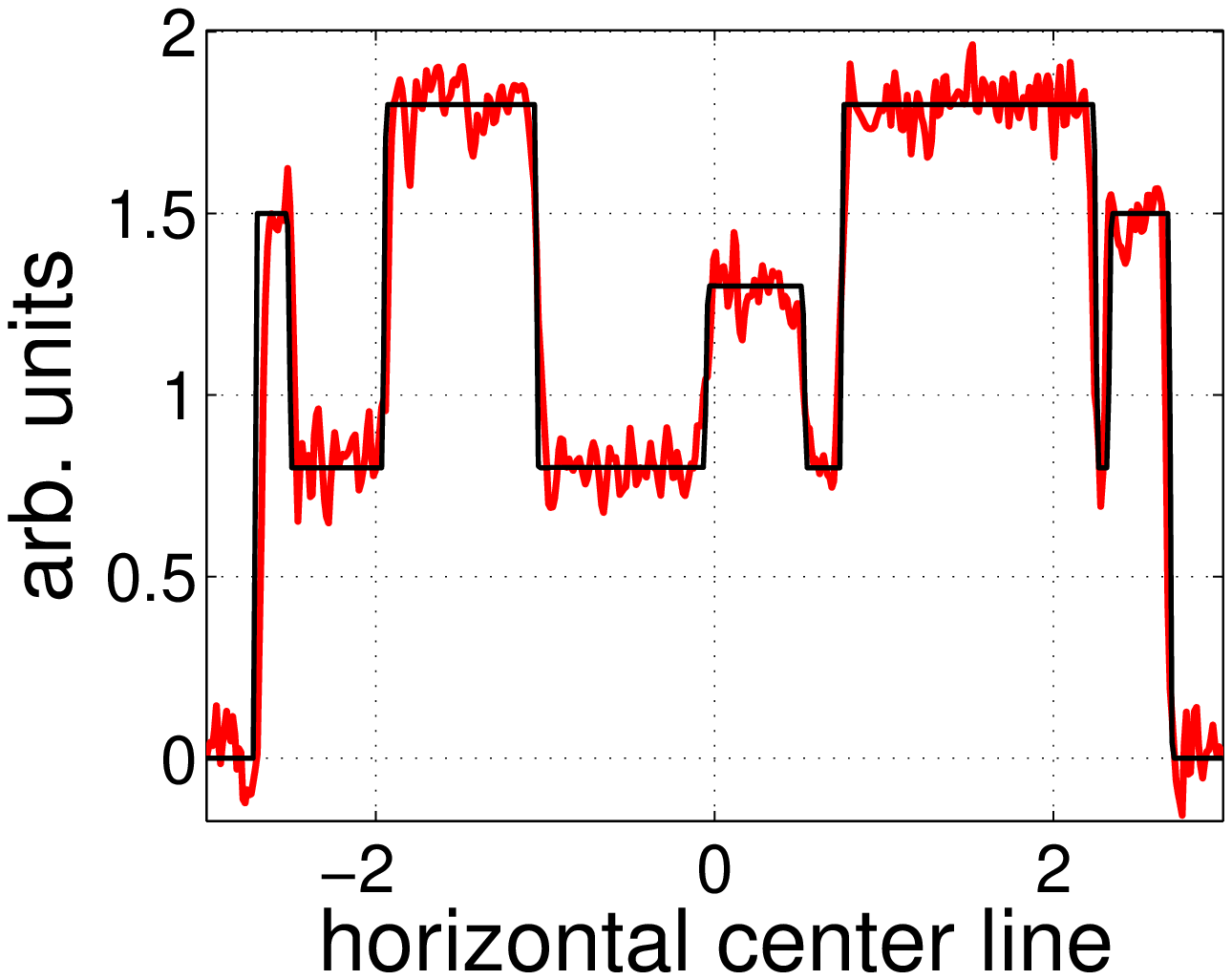}\\ 
\end{center}
\caption{\textbf{Numerical reconstruction using \req{wavefinite} 
from trace $\W f$ of the solution of the wave equation with $10\%$ noise
added.}}  
\label{fig:wave}
\end{figure}

In the following we present numerical results of our FBP algorithms for
reconstruction the phantom shown in the left picture in Figure
\ref{fig:phantom}, consisting of a superposition of characteristic functions
and one Gaussian kernel.  We calculated the data $\M f$ via numerical
integration and the operator $\W f = \di_\rT \Po f$ using
(\ref{ndsol}). Subsequently we added $5\%$ uniformly distributed noise to
$\M f$ and $10\%$ uniformly distributed noise to $\W f$.  The results for $N
= N_\ph = N_{\rM} = 300$ using the algorithms based on \req{m-lap},
\req{m-inv}, \req{m-inv-hilb} and \req{wavefinite} are depicted in Figures
\ref{fig:log}, \ref{fig:laplog} and \ref{fig:wave}.  All implementations
show good results although no explicit regularization strategy is
incorporated in order to the regularize the involved (mildly) ill posed
numerical differentiation.  In particular, \req{m-inv-hilb} and
\req{wavefinite} appear to be most insensitive to noise.  However, for noisy
data, the accuracy of FBP algorithms can be further improved by
incorporating a regularizing strategy similar to that used in \cite{HSS05}.
The derived identities in this article provide the mathematical foundation
for further development of FBP algorithms for the inversion from spherical
means and the inversion of the wave equation.

\bibliographystyle{amsplain}
\bibliography{twoD-1}
\end{document}